\documentclass[12pt]{article}
\usepackage{amssymb,amsfonts,amsmath, psfrag,eepic,colordvi,graphicx,epsfig,ytableau}
\usepackage{amssymb,latexsym,graphics,array}
\usepackage{MnSymbol}
\usepackage[enableskew]{youngtab}
\usepackage{float}
\usepackage{tikz,ifthen}
\usepackage{pgfplots}
\usetikzlibrary{math}
\topskip  
\newcommand{\rmnum}[1]{\romannumeral #1}

\numberwithin{equation}{section}
\newtheorem{theorem}{Theorem}[section]

\newtheorem{conjecture}[theorem]{Conjecture}

\newtheorem{lemma}[theorem]{Lemma}

\newtheorem{observation}[theorem]{Observation}

\begin{document}
	\parskip 6pt
	
	\pagenumbering{arabic}
	\def\sof{\hfill\rule{2mm}{2mm}}
	\def\ls{\leq}
	\def\gs{\geq}
	\def\SS{\mathcal S}
	\def\qq{{\bold q}}
	\def\MM{\mathcal M}
	\def\TT{\mathcal T}
	\def\EE{\mathcal E}
	\def\lsp{\mbox{lsp}}
	\def\rsp{\mbox{rsp}}
	\def\pf{\noindent {\it Proof.} }
	\def\mp{\mbox{pyramid}}
	\def\mb{\mbox{block}}
	\def\mc{\mbox{cross}}
	\def\qed{\hfill \rule{4pt}{7pt}}
	\def\block{\hfill \rule{5pt}{5pt}}
    \def\lr#1{\multicolumn{1}{|@{\hspace{.6ex}}c@{\hspace{.6ex}}|}{\raisebox{-.3ex}{$#1$}}}

	\begin{center}
		{\Large \bf    Proof of   a  conjecture on the  shape-Wilf-equivalence for  partially ordered patterns }
	\end{center}
	
	\begin{center}
	{\small Lintong Wang,  Sherry H.F. Yan$^{*}$\footnote{$^*$Corresponding author.}  \footnote{{\em E-mail address:}  hfy@zjnu.cn. }}

	 Department of Mathematics,
	Zhejiang Normal University\\
	Jinhua 321004, P.R. China		
\end{center}
	
\noindent {\bf Abstract.}   A   partially ordered pattern (abbreviated POP)   is a partially ordered set (poset) that
generalizes the notion of a pattern when we are not concerned with the relative order of some of its letters. The notion of  partially ordered patterns provides  a convenient language to deal with large  sets of permutation patterns. In analogy to the shape-Wilf-equivalence  for permutation patterns, Burstein-Han-Kitaev-Zhang initiated the study of the shape-Wilf-equivalence  for POPs which would   result  in the shape-Wilf-equivalence for large sets of permutation patterns.    The main objective of this paper is to confirm a recent   intriguing conjecture posed by Burstein-Han-Kitaev-Zhang concerning the shape-Wilf-equivalence for   POPs of length $k$. This is accomplished by establishing a   bijection  between two sets of  pattern-avoiding  transversals of a given Young diagram.  
 
\noindent {\bf Keywords}: pattern avoidance,  partially ordered pattern, shape-Wilf-equivalence.

\noindent {\bf AMS  Subject Classifications}: 05A05, 05C30

	
\section{Introduction}

 Let $\mathcal{S}_n$ denote the set of permutations of $[n]=\{1,2,\ldots, n\}$.
 Given a permutation $\pi=\pi_1\pi_2\ldots\pi_n \in \mathcal{S}_n$ and a permutation $\sigma=\sigma_1\sigma_2\ldots \sigma_k \in \mathcal{S}_k$,
an {\em occurrence} of $\sigma$ in $\pi$ is a subsequence $\pi_{i_1}\pi_{i_2}\cdots \pi_{i_k}$
of $\pi$ that is order isomorphic to $\sigma$.
We say $\pi$   {\em contains}  the pattern $\sigma$ if $\pi$ contains an occurrence of $\sigma$.
Otherwise, we say $\pi$ {\em avoids} the pattern $\sigma$  and
$\pi$ is {\em $\sigma$-avoiding}.  Given a set $P$ of permutation patterns,  a permutation $\pi$ is said to avoid  $P$ if $\pi$ avoids each pattern of $P$. Let $\mathcal{S}_n(P)$ denote the set of permutations in $\mathcal{S}_n$ that avoid $P$.  For two sets $P$ and $Q$ of permutation patterns, We say that $P$ is {\em Wilf-equivalent} to $Q$, denoted by $P\sim Q$, if the equality  $|\mathcal{S}_n(P)|=|\mathcal{S}_n(Q)|$  holds for all $n\geq 1$.

Pattern avoiding permutations were introduced by Knuth \cite{Knuth} in 1970 and first
systematically studied by Simion--Schmidt \cite{Simion}.
The   notion of pattern  avoidance in permutations  has been  extended to transversals of a Young diagram in  \cite{BWX}. Since then, 
  various results have been obtained for
pattern avoiding  tansversals  of  a Young diagram,  
	 see (\cite{Bousquet, Chan, Stankova, Yan2013,Yan2023, Zhou-Yan, Zhou-Yan-aoco}) and references therein.    

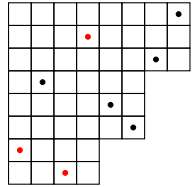
\begin{figure} [H]
\begin{center}
	\begin{tikzpicture}[scale = 0.3]
		\draw (0,0)rectangle(1,1);
		\draw (1,0)rectangle(2,1);
		\draw (2,0)rectangle(3,1);
		\draw (3,0)rectangle(4,1);
		\draw (0,1)rectangle(1,2);
		\draw (1,1)rectangle(2,2);
		\draw (2,1)rectangle(3,2);
		\draw (3,1)rectangle(4,2);
		\draw (0,2)rectangle(1,3);
		\draw (1,2)rectangle(2,3);
		\draw (2,2)rectangle(3,3);
		\draw (3,2)rectangle(4,3);
		\draw (4,2)rectangle(5,3);
		\draw (5,2)rectangle(6,3);
		\draw (0,3)rectangle(1,4);
		\draw (1,3)rectangle(2,4);
		\draw (2,3)rectangle(3,4);
		\draw (3,3)rectangle(4,4);
		\draw (4,3)rectangle(5,4);
		\draw (5,3)rectangle(6,4);
		\draw (0,4)rectangle(1,5);
		\draw (1,4)rectangle(2,5);
		\draw (2,4)rectangle(3,5);
		\draw (3,4)rectangle(4,5);
		\draw (4,4)rectangle(5,5);
		\draw (5,4)rectangle(6,5);
		\draw (0,5)rectangle(1,6);
		\draw (1,5)rectangle(2,6);
		\draw (2,5)rectangle(3,6);
		\draw (3,5)rectangle(4,6);
		\draw (4,5)rectangle(5,6);
		\draw (5,5)rectangle(6,6);
		\draw (6,5)rectangle(7,6);
		\draw (7,5)rectangle(8,6);
		\draw (0,6)rectangle(1,7);
		\draw (1,6)rectangle(2,7);
		\draw (2,6)rectangle(3,7);
		\draw (3,6)rectangle(4,7);
		\draw (4,6)rectangle(5,7);
		\draw (5,6)rectangle(6,7);
		\draw (6,6)rectangle(7,7);
		\draw (7,6)rectangle(8,7);
		\draw (0,7)rectangle(1,8);
		\draw (1,7)rectangle(2,8);
		\draw (2,7)rectangle(3,8);
		\draw (3,7)rectangle(4,8);
		\draw (4,7)rectangle(5,8);
		\draw (5,7)rectangle(6,8);
		\draw (6,7)rectangle(7,8);
		\draw (7,7)rectangle(8,8);
		\filldraw[red](0.5,1.5)circle(3pt);
		\filldraw[black](1.5,4.5)circle(3pt);
		\filldraw[red](2.5,0.5)circle(3pt);
		\filldraw[red](3.5,6.5)circle(3pt);
		\filldraw[black](4.5,3.5)circle(3pt);
		\filldraw[black](5.5,2.5)circle(3pt);
		\filldraw[black](6.5,5.5)circle(3pt);
		\filldraw[black](7.5,7.5)circle(3pt);
	\end{tikzpicture}
\end{center}
 \vspace{-0.5cm} 
	\caption{A transversal $T$ of the Young diagram $\lambda=(8,8,8,6,6,6,4,4)$. }\label{fig:tran}
 
\end{figure}

 In the following, let us first review some terminology related to   pattern avoiding  transversals. Throughout the  paper, we draw Young
 diagrams in English notation,  and   number columns from left to right and  rows from top to bottom.   We use the    square $(i,j)$ to represent  the square located in  row $i$   and column $j$. 
 A {\em  transversal} of  a Young diagram $\lambda=(\lambda_1, \lambda_2,   \ldots,  \lambda_n)$  with $\lambda_1\geq \lambda_2\geq \cdots \geq \lambda_n>0$ is a $01$-filling  of the squares of $\lambda$ with $1's$ and $0's$ such that every row and every column contains  exactly one $1$.   For better visibility, we will represent
 a $01$-filling by replacing a $1$ with a $\bullet$ but suppressing all occurrences of $0's$.   A transversal $T$ of the Young diagram $\lambda=(8,8,8,6,6,6,4,4)$  is displayed in  Figure \ref{fig:tran}.  
A permutation $\pi = \pi_1\pi_2\cdots\pi_n$ can be regarded as a transversal of the $n$ by $n$ square diagram,
in which the square $(\pi_i, i)$ is filled with a $1$ for all $1\leq i\leq n$
and all the other squares  are filled with $0's$. 
The transversal corresponding to the permutation $\pi$ is also called the
{\em permutation matrix} of $\pi$.
Let $\mathcal{S}_{\lambda}$ denote the set of transversals of the Young diagram $\lambda$.

Given a permutation $\alpha$ of $\mathcal{S}_m$,  let   $M$ be its permutation matrix. A transversal $T$ of a Young diagram $\lambda$ with $n$ columns and $n$ rows  will be said to contain  the pattern  $\alpha$ if there exists two subsets of the index set $[n]$, namely,  $R=\{r_1, r_2, \ldots, r_m\}$ and $C=\{c_1,  c_2,  \ldots,  c_m\}$, such that  the   matrix $M'$ on  the set $R$ of rows and the set  $C$ of columns is a copy of $M$ and each of the squares $(r_i,c_j)$ falls within the Young diagram.  In this context, we say that the matrix  $M'$ is isomorphic to   $\alpha$ or the matrix $M'$ is an occurrence of $\alpha$. 
Otherwise, we say that  $T$ {\em avoids}  the pattern $\alpha$ and $T$ is {\em $\alpha$-avoiding}.  For example, the submatrix on the $R=\{2,7,8\}$ and $C=\{1,3,4\}$ in the transversal  $T$  of  Figure \ref{fig:tran}  is an occurrence of pattern $231$. 
Given a set $P$ of patterns,  
let $\mathcal{S}_{\lambda}(P)$  denote the set of
transversals of the Young diagram $\lambda$ that avoid each pattern in $P$.  
 Given two sets $P$ and $Q$ of patterns, we say  that $P$ and $Q$  are {\em shape-Wilf-equivalent} if
  $|\mathcal{S}_{\lambda}(P)| = |\mathcal{S}_{\lambda}(Q)|$ holds  for any Young diagram  $\lambda$.   In this context, we write $P\sim_{s}Q$. If $\{\sigma\} \sim_{s} \{\tau\}$, we simply write $ \sigma \sim_s  \tau$.  Clearly, the shape-Wilf-equivalence would imply  the Wilf-equivalence. 

Given two permutations $\pi = \pi_1\pi_2\cdots \pi_n \in \mathcal{S}_n$
and $\sigma = \sigma_1\sigma_2\cdots \sigma_m \in \mathcal{S}_m$,
the {\em direct sum} of $\pi$ and $\sigma$, denoted by $\pi \oplus \sigma$, is the
permutation $\pi_1\pi_2\cdots \pi_n (\sigma_1+n)(\sigma_2+n)\cdots (\sigma_m+n)$.
 Backelin--West--Xin \cite{BWX}   obtained the following results concerning the shape-Wilf-equivalence  for permutation patterns.
\begin{theorem}\label{BWX1} (\cite{BWX}, Proposition  2.2)\label{BWX1}
	 For all $k\geq 1$, $I_k\sim_{s}J_k\sim_{s}  F_k$, where  $I_k = 12\cdots k$, $J_k = k\cdots21$ and $F_k=(k-1)\ldots 21 k $.
\end{theorem}
\begin{theorem}\label{BWX2} (\cite{BWX}, Proposition  2.3)
	For  any  patterns $\alpha$, $\beta$ and $\tau$, if $\alpha\sim_{s}\beta$, then  $\alpha\oplus \tau \sim_{s} \beta\oplus \tau$.
\end{theorem}
It should be mentioned that the cases $k=2$ and $k=3$ of Theorem \ref{BWX1}  have been   proved by  West \cite{West} and  Babson--West $\cite{BW}$, respectively.

Given  a set  $P=\{\tau^{(1)}, \tau^{(2)}, \ldots, \tau^{(m)}\}$ of  patterns  and a set $Q=\{\sigma^{(1)}, \sigma^{(2)}, \ldots, \sigma^{(n)}\}$ of  patterns,   the {\em direct sum} of $P$ and $Q$, denoted by $P\oplus Q$, is given by the set
$\{ \tau^{(i)}\oplus \sigma^{(j)}\mid 1\leq i\leq m, 1\leq j\leq n \}$.  Burstein-Han-Kitaev-Zhang \cite{Burstein} derived the following generalization of  
Theorem  \ref{BWX2}.
\begin{theorem}\label{Burstein} (\cite{Burstein}, Theorem 2)
	For  any sets $P$, $Q$ and $R$ of  patterns, if $P\sim_{s}Q$, then  $P\oplus R \sim_{s} Q\oplus R$.
\end{theorem}

  A {\em partially ordered pattern} (abbreviated POP) is a partially ordered set (poset) that
 generalizes the notion of a pattern when we are not concerned with the relative order of
 some of its letters, and therefore may represent a set of  patterns.   To be more specific, a POP $p$ of length $k$ is
 a poset with $k$ elements labeled by  $1, 2, \ldots, k$. Given a permutation $\pi\in \mathcal{S}_n$, we say that $\pi$ {\em contains} the POP $p$ if there exists a  
subsequence  $\pi_{i_1}\pi_{i_2}\cdots \pi_{i_k}$
of $\pi$ such that  that $s<_{p} t$ implies that $\pi_{i_s}<\pi_{i_t}$ for all $1\leq s, t\leq k$.   Otherwise, it is said to {\em avoid} the POP $p$ or to be  {\em $p$-avoiding}.  For example,  the POP $p$ =
\begin{tikzpicture}[baseline=(current bounding box.center), scale=0.5]
	 	\filldraw[black](0,0)circle(4pt);
	\node at(-0.5,0) {3};
	\filldraw[black](0,1.5)circle(4pt);
	\node at(-0.5,1.5) {1};
	\filldraw[black](1.5,0)circle(4pt);
	\node at(2,0) {2};
	\draw[black,very thick] (0,0)--(0,1.5);
\end{tikzpicture} represents all the patterns of length $3$ whose third element is smaller than its first element.   To be more specific,  the POP $p$ =
\begin{tikzpicture}[baseline=(current bounding box.center), scale=0.5]
\filldraw[black](0,0)circle(4pt);
\node at(-0.5,0) {3};
\filldraw[black](0,1.5)circle(4pt);
\node at(-0.5,1.5) {1};
\filldraw[black](1.5,0)circle(4pt);
\node at(2,0) {2};
\draw[black,very thick] (0,0)--(0,1.5);
\end{tikzpicture}
represents the patterns $231, 321, 312$ and the permutation $\pi$ avoids the POP $p$ =
\begin{tikzpicture}[baseline=(current bounding box.center), scale=0.5]
	\filldraw[black](0,0)circle(4pt);
	\node at(-0.5,0) {3};
	\filldraw[black](0,1.5)circle(4pt);
	\node at(-0.5,1.5) {1};
	\filldraw[black](1.5,0)circle(4pt);
	\node at(2,0) {2};
	\draw[black,very thick] (0,0)--(0,1.5);
\end{tikzpicture} if and only if it avoids the patterns $231, 321, 312$ simultaneously.
 The notion of  partially ordered patterns provides  a convenient language to deal with larger sets of permutation patterns.   Partially ordered patterns have been extensively exploited   in permutations and words, see \cite{Burstein, Gao, Han, Kitaev2005, Kitaev2007, Kitaev2003, Yap} and references therein.

 In analogy to the shape-Wilf-equivalence  for ordinary patterns,  Burstein-Han-Kitaev-Zhang \cite{Burstein}  initiated the study of the shape-Wilf-equivalence  for POPs. Given two POPs $p_1$ and $p_2$, let $P$ and $Q$ be the set of patterns that are represented by $p_1$ and $p_2$,  respectively. Then we say that $p_1$ is shape-Wilf-equivalent to $p_2$, denoted by $p_1\sim_s p_2$,  if $P\sim_s Q$.  Relying on Theorem \ref{Burstein}, the  shape-Wilf-equivalence  for POPs  allows  us to obtain  the shape-Wilf-equivalence  for  large sets of permutation patterns.
    In \cite{Burstein}, Burstein-Han-Kitaev-Zhang   derived  the following results on the  shape-Wilf-equivalence  for several  classes of POPs.   
 
 \begin{theorem}(\cite{Burstein}, Theorem 3)\label{thB1}
 For $k\geq 1$, let  $ x_1x_2 \ldots   x_k $ and $y_1y_2\ldots y_k$ be any two permutations of $[k]$. Then we have   
  \begin{tikzpicture}[baseline=(current bounding box.center), scale=0.5]
 	\filldraw[black](0,0)circle(4pt);
 	\node at(0,0.5) {$x_1$};
 	\filldraw[black](-2,-1.5)circle(4pt);
 	\node at(-2,-2) {$x_2$};
 	\filldraw[black](-1,-1.5)circle(4pt);
 	\node at(-1,-2) {$x_3$};
 	\filldraw[black](0.1,-1.5)circle(1.5pt);
 	\filldraw[black](0.5,-1.5)circle(1.5pt);
 	\filldraw[black](0.9,-1.5)circle(1.5pt);
 	\filldraw[black](2,-1.5)circle(4pt);
 	\node at(2,-2) {$x_{k}$};
 	\draw[black,very thick] (0,0)--(-2,-1.5);
 	\draw[black,very thick] (0,0)--(-1,-1.5);
 	\draw[black,very thick] (0,0)--(2,-1.5);
 \end{tikzpicture} 
 $\sim_s$
  \begin{tikzpicture}[baseline=(current bounding box.center), scale=0.5]
 	\filldraw[black](0,0)circle(4pt);
 	\node at(0,0.5) {$y_1$};
 	\filldraw[black](-2,-1.5)circle(4pt);
 	\node at(-2,-2) {$y_2$};
 	\filldraw[black](-1,-1.5)circle(4pt);
 	\node at(-1,-2) {$y_3$};
 	\filldraw[black](0.1,-1.5)circle(1.5pt);
 	\filldraw[black](0.5,-1.5)circle(1.5pt);
 	\filldraw[black](0.9,-1.5)circle(1.5pt);
 	\filldraw[black](2,-1.5)circle(4pt);
 	\node at(2,-2) {$y_{k}$};
 	\draw[black,very thick] (0,0)--(-2,-1.5);
 	\draw[black,very thick] (0,0)--(-1,-1.5);
 	\draw[black,very thick] (0,0)--(2,-1.5);
 \end{tikzpicture}.
 \end{theorem}
 
 \begin{theorem}(\cite{Burstein}, Theorem 4)\label{thB2}
 For $k\geq 1$,  let  $ x_1x_2 \ldots   x_k $ be any  permutation  of $[k]$. Then we have 
 
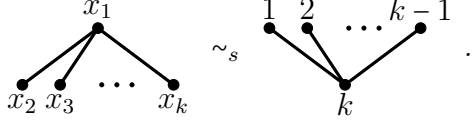
\begin{tikzpicture}[baseline=(current bounding box.center), scale=0.5]
 	\filldraw[black](0,0)circle(4pt);
 	\node at(0,0.5) {$x_1$};
 	\filldraw[black](-2,-1.5)circle(4pt);
 	\node at(-2,-2) {$x_2$};
 	\filldraw[black](-1,-1.5)circle(4pt);
 	\node at(-1,-2) {$x_3$};
 	\filldraw[black](0.1,-1.5)circle(1.5pt);
 	\filldraw[black](0.5,-1.5)circle(1.5pt);
 	\filldraw[black](0.9,-1.5)circle(1.5pt);
 	\filldraw[black](2,-1.5)circle(4pt);
 	\node at(2,-2) {$x_{k}$};
 	\draw[black,very thick] (0,0)--(-2,-1.5);
 	\draw[black,very thick] (0,0)--(-1,-1.5);
 	\draw[black,very thick] (0,0)--(2,-1.5);
 \end{tikzpicture} 
 $\sim_s$
\begin{tikzpicture}[baseline=(current bounding box.center), scale=0.5]
	\filldraw[black](0,0)circle(4pt);
	\node at(0,-0.5) {$k$};
	\filldraw[black](-2,1.5)circle(4pt);
	\node at(-2,2) {$1$};
	\filldraw[black](-1,1.5)circle(4pt);
	\node at(-1,2) {$2$};
	\filldraw[black](0.1,1.5)circle(1.5pt);
	\filldraw[black](0.5,1.5)circle(1.5pt);
	\filldraw[black](0.9,1.5)circle(1.5pt);
	\filldraw[black](2,1.5)circle(4pt);
	\node at(2,2) {$k-1$};
	\draw[black,very thick] (0,0)--(-2,1.5);
	\draw[black,very thick] (0,0)--(-1,1.5);
	\draw[black,very thick] (0,0)--(2,1.5);
\end{tikzpicture}. 
 \end{theorem}
  Burstein-Han-Kitaev-Zhang \cite{Burstein} further proposed the following intriguing     conjecture. 
 \begin{conjecture}(\cite{Burstein}, Conjecture 10)\label{con1}
 	For $k\geq 2$, we have 
 	$$ 
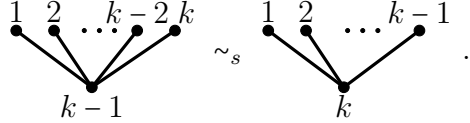
\begin{tikzpicture}[baseline=(current bounding box.center), scale=0.5]
 		\filldraw[black](0,0)circle(4pt);
 		\node at(0,-0.5) {$k-1$};
 		\filldraw[black](-2,1.5)circle(4pt);
 		\node at(-2,2) {$1$};
 		\filldraw[black](-1,1.5)circle(4pt);
 		\node at(-1,2) {$2$};
 		\filldraw[black](-0.2,1.5)circle(1.5pt);
 		\filldraw[black](0.2,1.5)circle(1.5pt);
 		\filldraw[black](0.6,1.5)circle(1.5pt);
 		\filldraw[black](1.2,1.5)circle(4pt);
 		\node at(1.2,2) {$k-2$};
 		\filldraw[black](2.2,1.5)circle(4pt);
 		\node at(2.5,2) {$k$};
 		\draw[black,very thick] (0,0)--(-2,1.5);
 		\draw[black,very thick] (0,0)--(-1,1.5);
 		\draw[black,very thick] (0,0)--(1.2,1.5);
 		\draw[black,very thick] (0,0)--(2.2,1.5);
 	\end{tikzpicture}
 	 \sim_s 
 	 \begin{tikzpicture}[baseline=(current bounding box.center), scale=0.5]
 	 	\filldraw[black](0,0)circle(4pt);
 	 	\node at(0,-0.5) {$k$};
 	 	\filldraw[black](-2,1.5)circle(4pt);
 	 	\node at(-2,2) {$1$};
 	 	\filldraw[black](-1,1.5)circle(4pt);
 	 	\node at(-1,2) {$2$};
 	 	\filldraw[black](0.1,1.5)circle(1.5pt);
 	 	\filldraw[black](0.5,1.5)circle(1.5pt);
 	 	\filldraw[black](0.9,1.5)circle(1.5pt);
 	 	\filldraw[black](2,1.5)circle(4pt);
 	 	\node at(2,2) {$k-1$};
 	 	\draw[black,very thick] (0,0)--(-2,1.5);
 	 	\draw[black,very thick] (0,0)--(-1,1.5);
 	 	\draw[black,very thick] (0,0)--(2,1.5);
 	 \end{tikzpicture}.
 	 $$
 \end{conjecture}
 
 It should be mentioned that the case $k=2$ of Conjecture \ref{con1} follows directly from Theorem \ref{BWX1} and the case $k=3$ of Conjecture \ref{con1} has been confirmed by   Burstein-Han- Kitaev-Zhang \cite{Burstein}.

 Let $P_k=$ \begin{tikzpicture}[baseline=(current bounding box.center), scale=0.5]
 	\filldraw[black](0,0)circle(4pt);
 	\node at(0,-0.5) {$k-1$};
 	\filldraw[black](-2,1.5)circle(4pt);
 	\node at(-2,2) {1};
 	\filldraw[black](-1,1.5)circle(4pt);
 	\node at(-1,2) {2};
 	\filldraw[black](-0.2,1.5)circle(1.5pt);
 	\filldraw[black](0.2,1.5)circle(1.5pt);
 	\filldraw[black](0.6,1.5)circle(1.5pt);
 	\filldraw[black](1.2,1.5)circle(4pt);
 	\node at(1.2,2) {$k-2$};
 	\filldraw[black](2.2,1.5)circle(4pt);
 	\node at(2.5,2) {$k$};
 	\draw[black,very thick] (0,0)--(-2,1.5);
 	\draw[black,very thick] (0,0)--(-1,1.5);
 	\draw[black,very thick] (0,0)--(1.2,1.5);
 	\draw[black,very thick] (0,0)--(2.2,1.5);
 \end{tikzpicture} and  $Q_k=$\begin{tikzpicture}[baseline=(current bounding box.center), scale=0.5]
 	\filldraw[black](0,0)circle(4pt);
 	\node at(0,0.5) {$k$};
 	\filldraw[black](-2,-1.5)circle(4pt);
 	\node at(-2,-2) {1};
 	\filldraw[black](-1,-1.5)circle(4pt);
 	\node at(-1,-2) {2};
 	\filldraw[black](0.1,-1.5)circle(1.5pt);
 	\filldraw[black](0.5,-1.5)circle(1.5pt);
 	\filldraw[black](0.9,-1.5)circle(1.5pt);
 	\filldraw[black](2,-1.5)circle(4pt);
 	\node at(2,-2) {$k-1$};
 	\draw[black,very thick] (0,0)--(-2,-1.5);
 	\draw[black,very thick] (0,0)--(-1,-1.5);
 	\draw[black,very thick] (0,0)--(2,-1.5);
 \end{tikzpicture}. 
 In view of Theorem   \ref{thB2}, in order to prove Conjecture \ref{con1}, it suffices to show that  $P_k\sim_s Q_k$.   The main objective of this   paper is to establish a bijection between $\mathcal{S}_\lambda(P_k)$ and $\mathcal{S}_\lambda(Q_k)$  for all $k\geq 3$ and for  any Young diagram $\lambda$, thereby confirming Conjecture \ref{con1}. 
 
\section{Proofs}
Let $k\geq 3$ and let $\lambda$ be a given Young diagram. 
In this section, we shall establish  a map   $\Phi: \mathcal{S}_{\lambda}(P_k)\rightarrow \mathcal{S}_{\lambda}(Q_k)$  that recursively transforms  every   occurrence    of $Q_k$ to  an occurrence  of $P_k$.  In order to show that the map $\Phi$ is a bijection,  we define  a map   $\Psi: \mathcal{S}_{\lambda}(Q_k)\rightarrow \mathcal{S}_{\lambda}(P_k)$  that recursively transforms  every  occurrence of $P_k$ to  an occurrence  of $Q_k$.  We shall show that the maps $\Phi$ and $\Psi$ are inverses of each other. 
\subsection{The map $\Phi$}
In this subsection, we first    define a transformation $\theta$  on the  submatrices of the transversals  which changes an occurrence  of $Q_k$ to  an occurrence  of $P_k$. The map $\Phi$ is then constructed relying on the transformation $\theta$.

{\noindent \bf The transformation $\theta$.}\\
Given a transversal $T\in \mathcal{S}_{\lambda}$,  
let $G$ be a  submatrix  in $T$ that is  isomorphic to  $Q_k$ at rows $r_1<r_2<\cdots <r_k$, and columns $c_1<c_2<\cdots<c_k$. Label the $1's$ of $G$ from left to right by $a_1, a_2, \ldots, a_k$, that is,  the $1$ positioned at   column $c_i$ is labeled by $a_i$ for all $1\leq i\leq k$.    Now we generate a  submatrix $\theta(G)$ isomorphic to $P_k$  with the same squares of $G$  by   distinguishing the following two cases. 
\begin{itemize}
	\item[Case \upshape(\rmnum{1}).]  If $a_{k-1}$ is  located at row $r_{k-1}$,  we cyclically shift the $1$ located at row $r_i$ to  row  $r_{i+2}$ for all $1\leq i\leq k$ with the convention that $r_{k+1}:=r_1$ and $r_{k+2}:=r_2$. Let $\theta(G)$ denote the resulting  submatrix.    See Figure \ref{fig:theta1}  for an illustration. 
	\item[Case \upshape(\rmnum{2}).] Otherwise, we cyclically shift the $1$ located at  row $r_i$ 	to   row $r_{i+1}$
	for all $1\leq i\leq k$ with the convention that $r_{k+1}:=r_1$.  Denote by $G'$  the resulting  submatrix. Let $\theta(G)$ denote the  submatrix isomorphic to $P_k$  obtained from $G'$ by swapping  the rightmost two columns  and leaving all the other squares fixed.  See Figure \ref{fig:theta2}  for an illustration. 
\end{itemize}
\begin{figure} [H]
\begin{center}
	\begin{tikzpicture}[scale = 0.4]
		\draw (0,0)rectangle(1,1);
		\draw (1,0)rectangle(2,1);
		\draw (2,0)rectangle(3,1);
		\draw (3,0)rectangle(4,1);
		\draw (0,1)rectangle(1,2);
		\draw (1,1)rectangle(2,2);
		\draw (2,1)rectangle(3,2);
		\draw (3,1)rectangle(4,2);
		\draw (0,2)rectangle(1,3);
		\draw (1,2)rectangle(2,3);
		\draw (2,2)rectangle(3,3);
		\draw (3,2)rectangle(4,3);
		\draw (4,2)rectangle(5,3);
		\draw (5,2)rectangle(6,3);
		\draw (0,3)rectangle(1,4);
		\draw (1,3)rectangle(2,4);
		\draw (2,3)rectangle(3,4);
		\draw (3,3)rectangle(4,4);
		\draw (4,3)rectangle(5,4);
		\draw (5,3)rectangle(6,4);
		\draw (0,4)rectangle(1,5);
		\draw (1,4)rectangle(2,5);
		\draw (2,4)rectangle(3,5);
		\draw (3,4)rectangle(4,5);
		\draw (4,4)rectangle(5,5);
		\draw (5,4)rectangle(6,5);
		\draw (0,5)rectangle(1,6);
		\draw (1,5)rectangle(2,6);
		\draw (2,5)rectangle(3,6);
		\draw (3,5)rectangle(4,6);
		\draw (4,5)rectangle(5,6);
		\draw (5,5)rectangle(6,6);
		\draw (6,5)rectangle(7,6);
		\draw (7,5)rectangle(8,6);
		\draw (0,6)rectangle(1,7);
		\draw (1,6)rectangle(2,7);
		\draw (2,6)rectangle(3,7);
		\draw (3,6)rectangle(4,7);
		\draw (4,6)rectangle(5,7);
		\draw (5,6)rectangle(6,7);
		\draw (6,6)rectangle(7,7);
		\draw (7,6)rectangle(8,7);
		\draw (0,7)rectangle(1,8);
		\draw (1,7)rectangle(2,8);
		\draw (2,7)rectangle(3,8);
		\draw (3,7)rectangle(4,8);
		\draw (4,7)rectangle(5,8);
		\draw (5,7)rectangle(6,8);
		\draw (6,7)rectangle(7,8);
		\draw (7,7)rectangle(8,8);
		\filldraw[red](0.5,4.5)circle(3pt);
		\filldraw[black](1.5,0.5)circle(3pt);
		\filldraw[red](2.5,7.5)circle(3pt);
		\filldraw[black](3.5,1.5)circle(3pt);
		\filldraw[red](4.5,3.5)circle(3pt);
		\filldraw[red](5.5,2.5)circle(3pt);
		\filldraw[black](6.5,6.5)circle(3pt);
		\filldraw[black](7.5,5.5)circle(3pt);
		\node[red] at(0.5,8.5){\footnotesize $c_{1}$};
		\node[red] at(2.5,8.5){\footnotesize $c_{2}$};
		\node[red] at(4.5,8.5){\footnotesize $c_{3}$};
		\node[red] at(5.5,8.5){\footnotesize $c_{4}$};
		\node[red] at(8.5,7.5){\footnotesize $r_{1}$};
		\node[red] at(6.5,4.5){\footnotesize $r_{2}$};
		\node[red] at(6.5,3.5){\footnotesize $r_{3}$};
		\node[red] at(6.5,2.5){\footnotesize $r_{4}$};
		
		\draw [->](9.3,4)--(11.5,4);
		
		\draw (13,0)rectangle(14,1);
		\draw (14,0)rectangle(15,1);
		\draw (15,0)rectangle(16,1);
		\draw (16,0)rectangle(17,1);
		\draw (13,1)rectangle(14,2);
		\draw (14,1)rectangle(15,2);
		\draw (15,1)rectangle(16,2);
		\draw (16,1)rectangle(17,2);
		\draw (13,2)rectangle(14,3);
		\draw (14,2)rectangle(15,3);
		\draw (15,2)rectangle(16,3);
		\draw (16,2)rectangle(17,3);
		\draw (17,2)rectangle(18,3);
		\draw (18,2)rectangle(19,3);
		\draw (13,3)rectangle(14,4);
		\draw (14,3)rectangle(15,4);
		\draw (15,3)rectangle(16,4);
		\draw (16,3)rectangle(17,4);
		\draw (17,3)rectangle(18,4);
		\draw (18,3)rectangle(19,4);
		\draw (13,4)rectangle(14,5);
		\draw (14,4)rectangle(15,5);
		\draw (15,4)rectangle(16,5);
		\draw (16,4)rectangle(17,5);
		\draw (17,4)rectangle(18,5);
		\draw (18,4)rectangle(19,5);
		\draw (13,5)rectangle(14,6);
		\draw (14,5)rectangle(15,6);
		\draw (15,5)rectangle(16,6);
		\draw (16,5)rectangle(17,6);
		\draw (17,5)rectangle(18,6);
		\draw (18,5)rectangle(19,6);
		\draw (19,5)rectangle(20,6);
		\draw (20,5)rectangle(21,6);
		\draw (13,6)rectangle(14,7);
		\draw (14,6)rectangle(15,7);
		\draw (15,6)rectangle(16,7);
		\draw (16,6)rectangle(17,7);
		\draw (17,6)rectangle(18,7);
		\draw (18,6)rectangle(19,7);
		\draw (19,6)rectangle(20,7);
		\draw (20,6)rectangle(21,7);
		\draw (13,7)rectangle(14,8);
		\draw (14,7)rectangle(15,8);
		\draw (15,7)rectangle(16,8);
		\draw (16,7)rectangle(17,8);
		\draw (17,7)rectangle(18,8);
		\draw (18,7)rectangle(19,8);
		\draw (19,7)rectangle(20,8);
		\draw (20,7)rectangle(21,8);
		\filldraw[red](13.5,2.5)circle(3pt);
		\filldraw[black](14.5,0.5)circle(3pt);
		\filldraw[red](15.5,3.5)circle(3pt);
		\filldraw[black](16.5,1.5)circle(3pt);
		\filldraw[red](17.5,7.5)circle(3pt);
		\filldraw[red](18.5,4.5)circle(3pt);
		\filldraw[black](19.5,6.5)circle(3pt);
		\filldraw[black](20.5,5.5)circle(3pt);
		\node[red] at(13.5,8.5){\footnotesize $c_{1}$};
		\node[red] at(15.5,8.5){\footnotesize $c_{2}$};
		\node[red] at(17.5,8.5){\footnotesize $c_{3}$};
		\node[red] at(18.5,8.5){\footnotesize $c_{4}$};
		\node[red] at(21.5,7.5){\footnotesize $r_{1}$};
		\node[red] at(19.5,4.5){\footnotesize $r_{2}$};
		\node[red] at(19.5,3.5){\footnotesize $r_{3}$};
		\node[red] at(19.5,2.5){\footnotesize $r_{4}$};
		\node[black] at(4,-2){\footnotesize $G$};
		\node[black] at(17,-2){\footnotesize $\theta(G)$};
		
	\end{tikzpicture}
\end{center}
 \vspace{-0.5cm} 
	\caption{An example of Case \upshape(\rmnum{1}) of the transformation $\theta$. }\label{fig:theta1}
\end{figure}

\begin{figure} [H]
\begin{center}
	\begin{tikzpicture}[scale = 0.4]
		\draw (0,0)rectangle(1,1);
		\draw (1,0)rectangle(2,1);
		\draw (2,0)rectangle(3,1);
		\draw (3,0)rectangle(4,1);
		\draw (0,1)rectangle(1,2);
		\draw (1,1)rectangle(2,2);
		\draw (2,1)rectangle(3,2);
		\draw (3,1)rectangle(4,2);
		\draw (0,2)rectangle(1,3);
		\draw (1,2)rectangle(2,3);
		\draw (2,2)rectangle(3,3);
		\draw (3,2)rectangle(4,3);
		\draw (4,2)rectangle(5,3);
		\draw (5,2)rectangle(6,3);
		\draw (0,3)rectangle(1,4);
		\draw (1,3)rectangle(2,4);
		\draw (2,3)rectangle(3,4);
		\draw (3,3)rectangle(4,4);
		\draw (4,3)rectangle(5,4);
		\draw (5,3)rectangle(6,4);
		\draw (0,4)rectangle(1,5);
		\draw (1,4)rectangle(2,5);
		\draw (2,4)rectangle(3,5);
		\draw (3,4)rectangle(4,5);
		\draw (4,4)rectangle(5,5);
		\draw (5,4)rectangle(6,5);
		\draw (0,5)rectangle(1,6);
		\draw (1,5)rectangle(2,6);
		\draw (2,5)rectangle(3,6);
		\draw (3,5)rectangle(4,6);
		\draw (4,5)rectangle(5,6);
		\draw (5,5)rectangle(6,6);
		\draw (6,5)rectangle(7,6);
		\draw (7,5)rectangle(8,6);
		\draw (0,6)rectangle(1,7);
		\draw (1,6)rectangle(2,7);
		\draw (2,6)rectangle(3,7);
		\draw (3,6)rectangle(4,7);
		\draw (4,6)rectangle(5,7);
		\draw (5,6)rectangle(6,7);
		\draw (6,6)rectangle(7,7);
		\draw (7,6)rectangle(8,7);
		\draw (0,7)rectangle(1,8);
		\draw (1,7)rectangle(2,8);
		\draw (2,7)rectangle(3,8);
		\draw (3,7)rectangle(4,8);
		\draw (4,7)rectangle(5,8);
		\draw (5,7)rectangle(6,8);
		\draw (6,7)rectangle(7,8);
		\draw (7,7)rectangle(8,8);
		\filldraw[red](0.5,3.5)circle(3pt);
		\filldraw[black](1.5,0.5)circle(3pt);
		\filldraw[red](2.5,7.5)circle(3pt);
		\filldraw[black](3.5,1.5)circle(3pt);
		\filldraw[red](4.5,4.5)circle(3pt);
		\filldraw[red](5.5,2.5)circle(3pt);
		\filldraw[black](6.5,6.5)circle(3pt);
		\filldraw[black](7.5,5.5)circle(3pt);
		\node[red] at(0.5,8.5){\footnotesize $c_{1}$};
		\node[red] at(2.5,8.5){\footnotesize $c_{2}$};
		\node[red] at(4.5,8.5){\footnotesize $c_{3}$};
		\node[red] at(5.5,8.5){\footnotesize $c_{4}$};
		\node[red] at(8.5,7.5){\footnotesize $r_{1}$};
		\node[red] at(6.5,4.5){\footnotesize $r_{2}$};
		\node[red] at(6.5,3.5){\footnotesize $r_{3}$};
		\node[red] at(6.5,2.5){\footnotesize $r_{4}$};
		
		\draw [->](9,4)--(11.5,4);
		
		\draw (13,0)rectangle(14,1);
		\draw (14,0)rectangle(15,1);
		\draw (15,0)rectangle(16,1);
		\draw (16,0)rectangle(17,1);
		\draw (13,1)rectangle(14,2);
		\draw (14,1)rectangle(15,2);
		\draw (15,1)rectangle(16,2);
		\draw (16,1)rectangle(17,2);
		\draw (13,2)rectangle(14,3);
		\draw (14,2)rectangle(15,3);
		\draw (15,2)rectangle(16,3);
		\draw (16,2)rectangle(17,3);
		\draw (17,2)rectangle(18,3);
		\draw (18,2)rectangle(19,3);
		\draw (13,3)rectangle(14,4);
		\draw (14,3)rectangle(15,4);
		\draw (15,3)rectangle(16,4);
		\draw (16,3)rectangle(17,4);
		\draw (17,3)rectangle(18,4);
		\draw (18,3)rectangle(19,4);
		\draw (13,4)rectangle(14,5);
		\draw (14,4)rectangle(15,5);
		\draw (15,4)rectangle(16,5);
		\draw (16,4)rectangle(17,5);
		\draw (17,4)rectangle(18,5);
		\draw (18,4)rectangle(19,5);
		\draw (13,5)rectangle(14,6);
		\draw (14,5)rectangle(15,6);
		\draw (15,5)rectangle(16,6);
		\draw (16,5)rectangle(17,6);
		\draw (17,5)rectangle(18,6);
		\draw (18,5)rectangle(19,6);
		\draw (19,5)rectangle(20,6);
		\draw (20,5)rectangle(21,6);
		\draw (13,6)rectangle(14,7);
		\draw (14,6)rectangle(15,7);
		\draw (15,6)rectangle(16,7);
		\draw (16,6)rectangle(17,7);
		\draw (17,6)rectangle(18,7);
		\draw (18,6)rectangle(19,7);
		\draw (19,6)rectangle(20,7);
		\draw (20,6)rectangle(21,7);
		\draw (13,7)rectangle(14,8);
		\draw (14,7)rectangle(15,8);
		\draw (15,7)rectangle(16,8);
		\draw (16,7)rectangle(17,8);
		\draw (17,7)rectangle(18,8);
		\draw (18,7)rectangle(19,8);
		\draw (19,7)rectangle(20,8);
		\draw (20,7)rectangle(21,8);
		\filldraw[red](13.5,2.5)circle(3pt);
		\filldraw[black](14.5,0.5)circle(3pt);
		\filldraw[red](15.5,4.5)circle(3pt);
		\filldraw[black](16.5,1.5)circle(3pt);
		\filldraw[red](17.5,3.5)circle(3pt);
		\filldraw[red](18.5,7.5)circle(3pt);
		\filldraw[black](19.5,6.5)circle(3pt);
		\filldraw[black](20.5,5.5)circle(3pt);
		\node[red] at(13.5,8.5){\footnotesize $c_{1}$};
		\node[red] at(15.5,8.5){\footnotesize $c_{2}$};
		\node[red] at(17.5,8.5){\footnotesize $c_{3}$};
		\node[red] at(18.5,8.5){\footnotesize $c_{4}$};
		\node[red] at(21.5,7.5){\footnotesize $r_{1}$};
		\node[red] at(19.5,4.5){\footnotesize $r_{2}$};
		\node[red] at(19.5,3.5){\footnotesize $r_{3}$};
		\node[red] at(19.5,2.5){\footnotesize $r_{4}$};
		
		\draw [->](22,4)--(24.5,4);
		
		\draw (26,0)rectangle(27,1);
		\draw (27,0)rectangle(28,1);
		\draw (28,0)rectangle(29,1);
		\draw (29,0)rectangle(30,1);
		\draw (26,1)rectangle(27,2);
		\draw (27,1)rectangle(28,2);
		\draw (28,1)rectangle(29,2);
		\draw (29,1)rectangle(30,2);
		\draw (26,2)rectangle(27,3);
		\draw (27,2)rectangle(28,3);
		\draw (28,2)rectangle(29,3);
		\draw (29,2)rectangle(30,3);
		\draw (30,2)rectangle(31,3);
		\draw (31,2)rectangle(32,3);
		\draw (26,3)rectangle(27,4);
		\draw (27,3)rectangle(28,4);
		\draw (28,3)rectangle(29,4);
		\draw (29,3)rectangle(30,4);
		\draw (30,3)rectangle(31,4);
		\draw (31,3)rectangle(32,4);
		\draw (26,4)rectangle(27,5);
		\draw (27,4)rectangle(28,5);
		\draw (28,4)rectangle(29,5);
		\draw (29,4)rectangle(30,5);
		\draw (30,4)rectangle(31,5);
		\draw (31,4)rectangle(32,5);
		\draw (26,5)rectangle(27,6);
		\draw (27,5)rectangle(28,6);
		\draw (28,5)rectangle(29,6);
		\draw (29,5)rectangle(30,6);
		\draw (30,5)rectangle(31,6);
		\draw (31,5)rectangle(32,6);
		\draw (32,5)rectangle(33,6);
		\draw (33,5)rectangle(34,6);
		\draw (26,6)rectangle(27,7);
		\draw (27,6)rectangle(28,7);
		\draw (28,6)rectangle(29,7);
		\draw (29,6)rectangle(30,7);
		\draw (30,6)rectangle(31,7);
		\draw (31,6)rectangle(32,7);
		\draw (32,6)rectangle(33,7);
		\draw (33,6)rectangle(34,7);
		\draw (26,7)rectangle(27,8);
		\draw (27,7)rectangle(28,8);
		\draw (28,7)rectangle(29,8);
		\draw (29,7)rectangle(30,8);
		\draw (30,7)rectangle(31,8);
		\draw (31,7)rectangle(32,8);
		\draw (32,7)rectangle(33,8);
		\draw (33,7)rectangle(34,8);
		\filldraw[red](26.5,2.5)circle(3pt);
		\filldraw[black](27.5,0.5)circle(3pt);
		\filldraw[red](28.5,4.5)circle(3pt);
		\filldraw[black](29.5,1.5)circle(3pt);
		\filldraw[red](30.5,7.5)circle(3pt);
		\filldraw[red](31.5,3.5)circle(3pt);
		\filldraw[black](32.5,6.5)circle(3pt);
		\filldraw[black](33.5,5.5)circle(3pt);
		\node[red] at(26.5,8.5){\footnotesize $c_{1}$};
		\node[red] at(28.5,8.5){\footnotesize $c_{2}$};
		\node[red] at(30.5,8.5){\footnotesize $c_{3}$};
		\node[red] at(31.5,8.5){\footnotesize $c_{4}$};
		\node[red] at(34.5,7.5){\footnotesize $r_{1}$};
		\node[red] at(32.5,4.5){\footnotesize $r_{2}$};
		\node[red] at(32.5,3.5){\footnotesize $r_{3}$};
		\node[red] at(32.5,2.5){\footnotesize $r_{4}$};
		\node[black] at(4,-2){\footnotesize $G$};
		\node[black] at(17,-2){\footnotesize $G'$};
		\node[black] at(30,-2){\footnotesize $\theta(G)$};

	\end{tikzpicture}
\end{center}
 \vspace{-0.8cm} 
	\caption{An example of Case \upshape(\rmnum{2}) of the transformation $\theta$. }\label{fig:theta2}
\end{figure}

{\noindent \bf The map  $\Phi: \mathcal{S}_{\lambda}(P_k)\rightarrow \mathcal{S}_{\lambda}(Q_k)$ .}\\
Suppose that   $T$ is a transversal  in  $ \mathcal{S}_{\lambda}(P_k)$. 
 \begin{itemize}
	\item[Step 1.] If  $T$ contains no $Q_k$, end.
	\item[Step 2.]  Find the lowest square $b_1$ containing $1$, such that there is a $Q_k$ in $T$ in which $b_1$ is the lowest  $1$. 
	
		\item[Step 3.]  Find the lowest squares $b_2, b_3, \ldots, b_k$ (listed from bottom  to top) containing   $1's$ that are located above and  to the left of $b_1$.
		  Then we get a submatrix  $G$ isomorphic to $Q_k$. 
			\item[Step 4.]  Replace  the submatrix $G$  with $\theta(G)$ and leave all the other squares fixed. 	 
	\item[Step 5.]  Repeat the above procedure  until  there is no $Q_k$ left. Then we get a $\Phi(T)\in  \mathcal{S}_{\lambda}(Q_k)$.
\end{itemize}

For our convenience,    Steps 2--4 of the algorithm $\Phi$  is called  {\em the transformation $\phi$}.
Let  $t$ be a nonnegative  integer. Denote by $\phi^{t}(T)$ the resulting transversal  obtained  by the $t$-th application of $\phi $ to the $T$. Clearly, we have $\phi^{0}(T)=T$.

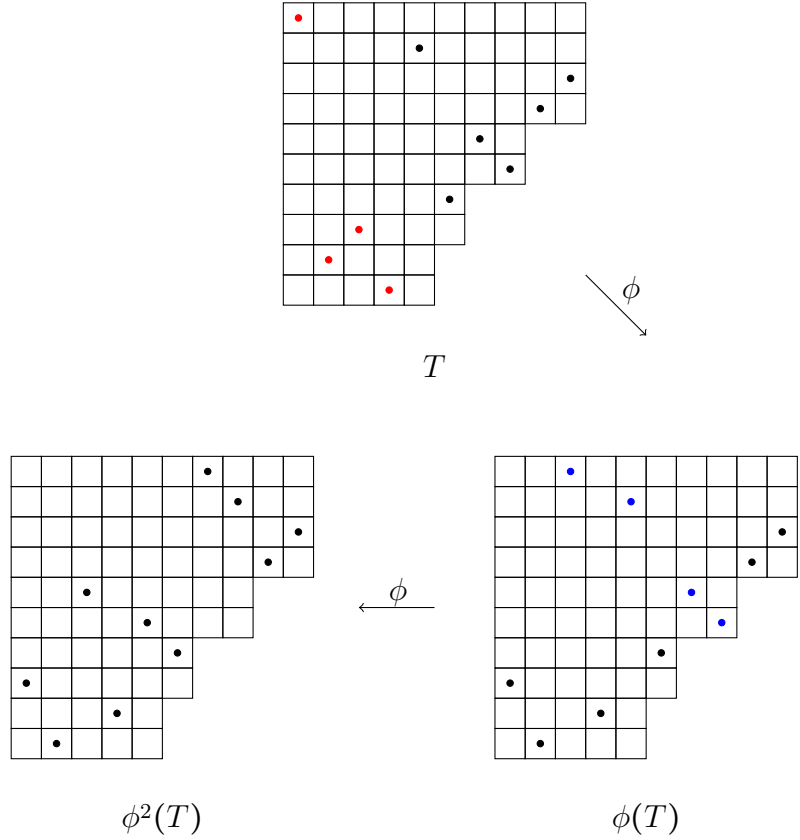
\begin{figure} [H]
\begin{center}
	\begin{tikzpicture}[scale = 0.4]
		\draw (0,0)rectangle(1,1);
		\draw (1,0)rectangle(2,1);
		\draw (2,0)rectangle(3,1);
		\draw (3,0)rectangle(4,1);
		\draw (4,0)rectangle(5,1);
		\draw (0,1)rectangle(1,2);
		\draw (1,1)rectangle(2,2);
		\draw (2,1)rectangle(3,2);
		\draw (3,1)rectangle(4,2);
		\draw (4,1)rectangle(5,2);
		\draw (0,2)rectangle(1,3);
		\draw (1,2)rectangle(2,3);
		\draw (2,2)rectangle(3,3);
		\draw (3,2)rectangle(4,3);
		\draw (4,2)rectangle(5,3);
		\draw (5,2)rectangle(6,3);
		\draw (0,3)rectangle(1,4);
		\draw (1,3)rectangle(2,4);
		\draw (2,3)rectangle(3,4);
		\draw (3,3)rectangle(4,4);
		\draw (4,3)rectangle(5,4);
		\draw (5,3)rectangle(6,4);
		\draw (0,4)rectangle(1,5);
		\draw (1,4)rectangle(2,5);
		\draw (2,4)rectangle(3,5);
		\draw (3,4)rectangle(4,5);
		\draw (4,4)rectangle(5,5);
		\draw (5,4)rectangle(6,5);
		\draw (6,4)rectangle(7,5);
		\draw (7,4)rectangle(8,5);
		\draw (0,5)rectangle(1,6);
		\draw (1,5)rectangle(2,6);
		\draw (2,5)rectangle(3,6);
		\draw (3,5)rectangle(4,6);
		\draw (4,5)rectangle(5,6);
		\draw (5,5)rectangle(6,6);
		\draw (6,5)rectangle(7,6);
		\draw (7,5)rectangle(8,6);
		\draw (0,6)rectangle(1,7);
		\draw (1,6)rectangle(2,7);
		\draw (2,6)rectangle(3,7);
		\draw (3,6)rectangle(4,7);
		\draw (4,6)rectangle(5,7);
		\draw (5,6)rectangle(6,7);
		\draw (6,6)rectangle(7,7);
		\draw (7,6)rectangle(8,7);
		\draw (8,6)rectangle(9,7);
		\draw (9,6)rectangle(10,7);
		\draw (0,7)rectangle(1,8);
		\draw (1,7)rectangle(2,8);
		\draw (2,7)rectangle(3,8);
		\draw (3,7)rectangle(4,8);
		\draw (4,7)rectangle(5,8);
		\draw (5,7)rectangle(6,8);
		\draw (6,7)rectangle(7,8);
		\draw (7,7)rectangle(8,8);
		\draw (8,7)rectangle(9,8);
		\draw (9,7)rectangle(10,8);
		\draw (0,8)rectangle(1,9);
		\draw (1,8)rectangle(2,9);
		\draw (2,8)rectangle(3,9);
		\draw (3,8)rectangle(4,9);
		\draw (4,8)rectangle(5,9);
		\draw (5,8)rectangle(6,9);
		\draw (6,8)rectangle(7,9);
		\draw (7,8)rectangle(8,9);
		\draw (8,8)rectangle(9,9);
		\draw (9,8)rectangle(10,9);
		\draw (0,9)rectangle(1,10);
		\draw (1,9)rectangle(2,10);
		\draw (2,9)rectangle(3,10);
		\draw (3,9)rectangle(4,10);
		\draw (4,9)rectangle(5,10);
		\draw (5,9)rectangle(6,10);
		\draw (6,9)rectangle(7,10);
		\draw (7,9)rectangle(8,10);
		\draw (8,9)rectangle(9,10);
		\draw (9,9)rectangle(10,10);
		\filldraw[black](0.5,2.5)circle(3pt);
		\filldraw[black](1.5,0.5)circle(3pt);
		\filldraw[black](2.5,5.5)circle(3pt);
		\filldraw[black](3.5,1.5)circle(3pt);
		\filldraw[black](4.5,4.5)circle(3pt);
		\filldraw[black](5.5,3.5)circle(3pt);
		\filldraw[black](6.5,9.5)circle(3pt);
		\filldraw[black](7.5,8.5)circle(3pt);
		\filldraw[black](8.5,6.5)circle(3pt);
		\filldraw[black](9.5,7.5)circle(3pt);
		\node at(5,-2){$\phi^{2}(T)$};
		
		\draw [<-](11.5,5)--(14,5);
		\node at(12.8,5.5){$\phi$};
		
		\draw (16,0)rectangle(17,1);
		\draw (17,0)rectangle(18,1);
		\draw (18,0)rectangle(19,1);
		\draw (19,0)rectangle(20,1);
		\draw (20,0)rectangle(21,1);
		\draw (16,1)rectangle(17,2);
		\draw (17,1)rectangle(18,2);
		\draw (18,1)rectangle(19,2);
		\draw (19,1)rectangle(20,2);
		\draw (20,1)rectangle(21,2);
		\draw (16,2)rectangle(17,3);
		\draw (17,2)rectangle(18,3);
		\draw (18,2)rectangle(19,3);
		\draw (19,2)rectangle(20,3);
		\draw (20,2)rectangle(21,3);
		\draw (21,2)rectangle(22,3);
		\draw (16,3)rectangle(17,4);
		\draw (17,3)rectangle(18,4);
		\draw (18,3)rectangle(19,4);
		\draw (19,3)rectangle(20,4);
		\draw (20,3)rectangle(21,4);
		\draw (21,3)rectangle(22,4);
		\draw (16,4)rectangle(17,5);
		\draw (17,4)rectangle(18,5);
		\draw (18,4)rectangle(19,5);
		\draw (19,4)rectangle(20,5);
		\draw (20,4)rectangle(21,5);
		\draw (21,4)rectangle(22,5);
		\draw (22,4)rectangle(23,5);
		\draw (23,4)rectangle(24,5);
		\draw (16,5)rectangle(17,6);
		\draw (17,5)rectangle(18,6);
		\draw (18,5)rectangle(19,6);
		\draw (19,5)rectangle(20,6);
		\draw (20,5)rectangle(21,6);
		\draw (21,5)rectangle(22,6);
		\draw (22,5)rectangle(23,6);
		\draw (23,5)rectangle(24,6);
		\draw (16,6)rectangle(17,7);
		\draw (17,6)rectangle(18,7);
		\draw (18,6)rectangle(19,7);
		\draw (19,6)rectangle(20,7);
		\draw (20,6)rectangle(21,7);
		\draw (21,6)rectangle(22,7);
		\draw (22,6)rectangle(23,7);
		\draw (23,6)rectangle(24,7);
		\draw (24,6)rectangle(25,7);
		\draw (25,6)rectangle(26,7);
		\draw (16,7)rectangle(17,8);
		\draw (17,7)rectangle(18,8);
		\draw (18,7)rectangle(19,8);
		\draw (19,7)rectangle(20,8);
		\draw (20,7)rectangle(21,8);
		\draw (21,7)rectangle(22,8);
		\draw (22,7)rectangle(23,8);
		\draw (23,7)rectangle(24,8);
		\draw (24,7)rectangle(25,8);
		\draw (25,7)rectangle(26,8);
		\draw (16,8)rectangle(17,9);
		\draw (17,8)rectangle(18,9);
		\draw (18,8)rectangle(19,9);
		\draw (19,8)rectangle(20,9);
		\draw (20,8)rectangle(21,9);
		\draw (21,8)rectangle(22,9);
		\draw (22,8)rectangle(23,9);
		\draw (23,8)rectangle(24,9);
		\draw (24,8)rectangle(25,9);
		\draw (25,8)rectangle(26,9);
		\draw (16,9)rectangle(17,10);
		\draw (17,9)rectangle(18,10);
		\draw (18,9)rectangle(19,10);
		\draw (19,9)rectangle(20,10);
		\draw (20,9)rectangle(21,10);
		\draw (21,9)rectangle(22,10);
		\draw (22,9)rectangle(23,10);
		\draw (23,9)rectangle(24,10);
		\draw (24,9)rectangle(25,10);
		\draw (25,9)rectangle(26,10);
		\filldraw[black](16.5,2.5)circle(3pt);
		\filldraw[black](17.5,0.5)circle(3pt);
		\filldraw[blue](18.5,9.5)circle(3pt);
		\filldraw[black](19.5,1.5)circle(3pt);
		\filldraw[blue](20.5,8.5)circle(3pt);
		\filldraw[black](21.5,3.5)circle(3pt);
		\filldraw[blue](22.5,5.5)circle(3pt);
		\filldraw[blue](23.5,4.5)circle(3pt);
		\filldraw[black](24.5,6.5)circle(3pt);
		\filldraw[black](25.5,7.5)circle(3pt);
		\node at(21,-2){$\phi(T)$};
		
		\draw [<-](21,14)--(19,16);
		\node at(20.5,15.5){$\phi$};
		
		\draw (9,15)rectangle(10,16);
		\draw (10,15)rectangle(11,16);
		\draw (11,15)rectangle(12,16);
		\draw (12,15)rectangle(13,16);
		\draw (13,15)rectangle(14,16);
		\draw (9,16)rectangle(10,17);
		\draw (10,16)rectangle(11,17);
		\draw (11,16)rectangle(12,17);
		\draw (12,16)rectangle(13,17);
		\draw (13,16)rectangle(14,17);
		\draw (9,17)rectangle(10,18);
		\draw (10,17)rectangle(11,18);
		\draw (11,17)rectangle(12,18);
		\draw (12,17)rectangle(13,18);
		\draw (13,17)rectangle(14,18);
		\draw (14,17)rectangle(15,18);
		\draw (9,18)rectangle(10,19);
		\draw (10,18)rectangle(11,19);
		\draw (11,18)rectangle(12,19);
		\draw (12,18)rectangle(13,19);
		\draw (13,18)rectangle(14,19);
		\draw (14,18)rectangle(15,19);
		\draw (9,19)rectangle(10,20);
		\draw (10,19)rectangle(11,20);
		\draw (11,19)rectangle(12,20);
		\draw (12,19)rectangle(13,20);
		\draw (13,19)rectangle(14,20);
		\draw (14,19)rectangle(15,20);
		\draw (15,19)rectangle(16,20);
		\draw (16,19)rectangle(17,20);
		\draw (9,20)rectangle(10,21);
		\draw (10,20)rectangle(11,21);
		\draw (11,20)rectangle(12,21);
		\draw (12,20)rectangle(13,21);
		\draw (13,20)rectangle(14,21);
		\draw (14,20)rectangle(15,21);
		\draw (15,20)rectangle(16,21);
		\draw (16,20)rectangle(17,21);
		\draw (9,21)rectangle(10,22);
		\draw (10,21)rectangle(11,22);
		\draw (11,21)rectangle(12,22);
		\draw (12,21)rectangle(13,22);
		\draw (13,21)rectangle(14,22);
		\draw (14,21)rectangle(15,22);
		\draw (15,21)rectangle(16,22);
		\draw (16,21)rectangle(17,22);
		\draw (17,21)rectangle(18,22);
		\draw (18,21)rectangle(19,22);
		\draw (9,22)rectangle(10,23);
		\draw (10,22)rectangle(11,23);
		\draw (11,22)rectangle(12,23);
		\draw (12,22)rectangle(13,23);
		\draw (13,22)rectangle(14,23);
		\draw (14,22)rectangle(15,23);
		\draw (15,22)rectangle(16,23);
		\draw (16,22)rectangle(17,23);
		\draw (17,22)rectangle(18,23);
		\draw (18,22)rectangle(19,23);
		\draw (9,23)rectangle(10,24);
		\draw (10,23)rectangle(11,24);
		\draw (11,23)rectangle(12,24);
		\draw (12,23)rectangle(13,24);
		\draw (13,23)rectangle(14,24);
		\draw (14,23)rectangle(15,24);
		\draw (15,23)rectangle(16,24);
		\draw (16,23)rectangle(17,24);
		\draw (17,23)rectangle(18,24);
		\draw (18,23)rectangle(19,24);
		\draw (9,24)rectangle(10,25);
		\draw (10,24)rectangle(11,25);
		\draw (11,24)rectangle(12,25);
		\draw (12,24)rectangle(13,25);
		\draw (13,24)rectangle(14,25);
		\draw (14,24)rectangle(15,25);
		\draw (15,24)rectangle(16,25);
		\draw (16,24)rectangle(17,25);
		\draw (17,24)rectangle(18,25);
		\draw (18,24)rectangle(19,25);
		\filldraw[red](9.5,24.5)circle(3pt);
		\filldraw[red](10.5,16.5)circle(3pt);
		\filldraw[red](11.5,17.5)circle(3pt);
		\filldraw[red](12.5,15.5)circle(3pt);
		\filldraw[black](13.5,23.5)circle(3pt);
		\filldraw[black](14.5,18.5)circle(3pt);
		\filldraw[black](15.5,20.5)circle(3pt);
		\filldraw[black](16.5,19.5)circle(3pt);
		\filldraw[black](17.5,21.5)circle(3pt);
		\filldraw[black](18.5,22.5)circle(3pt);
		\node at(14,13){$T$};
		
	\end{tikzpicture}
\end{center}
 \vspace{-0.5cm} 
 \caption{An example of the map $\Phi$. }\label{fig:Phi}
\end{figure}

For example, let $T$ be a transversal as shown in Figure \ref{fig:Phi}. When we apply the transformation $\phi$ to $T$, the selected submatrix $G$ isomorphic to $Q_4$ is the submatrix  whose $1's$ are colored by red as shown in Figure \ref{fig:Phi}. By replacing $G$ with $\theta(G)$, we get a transversal $\phi(T)$.  In the second application of $\phi$ to $T$,  the selected submatrix $G'$ isomorphic to $Q_4$ in $\phi(T)$ is the  submatrix  whose $1's$ are colored by blue as shown in Figure \ref{fig:Phi}. By replacing $G'$ with $\theta(G)$ in $\phi(T)$, we eventually get a transversal $\phi^2(T)$ containing no $Q_k$.

Now we proceed to prove that the transformation $\phi$ has the following celebrated properties. We begin with some necessary   definitions and notations.

Let $t\geq 1$ and $T\in  \mathcal{S}_{\lambda}(P_k)$. If $T$ contains $Q_k$,  
suppose that at the $t$-th application of $\phi$ to $T$, the selected  submatrix   isomorphic to $Q_k$ is given by $G$ in   which  the squares containing $1's$ of   $G$  are labeled  by $b_1, b_2, \ldots, b_k$ from bottom to top,  and  the $1's$ of $G$ are labeled  by $a_1, a_2, \ldots, a_k$ from left to right. Clearly, $b_1$ is the square containing the lowest $1$ of $G$ and $b_k$ is the square containing the topmost  $1$ of $G$.  Label the $1's$ of  the  submatrix $\theta(G)$ from left to right by $g_1, g_2, \ldots, g_k$.

\begin{figure}[H]
	\begin{center}
		\begin{tikzpicture}[scale = 0.4]
			\fill[lightgray](0,0)rectangle(13,4);
			\fill[lightgray](9,0)rectangle(13,10);
			\draw[black](0,0)--(13,0);
			\draw[black](13,0)--(13,10);
			\draw[black](0,10)--(13,10);
			\draw[black](9,10)--(9,4);
			\draw[black](9,4)--(0,4);
			\filldraw[black](13,0)circle(4pt);
			\filldraw[black](13,6)circle(4pt);
			\filldraw[black](6,10)circle(4pt);
			\filldraw[black](4,4)circle(4pt);
			\filldraw[black](9,8)circle(4pt);
			\filldraw[black](9,10)circle(4pt);
			\filldraw[black](3,0)circle(4pt);
			\node at(15,0) {$b_{1}(a_{k})$};
			\node at(14,6) {$g_{k}$};
			\node at(9,11) {$g_{k-1}$};
			\node at(3,-1) {$g_{s}$};
			\node at(6,11) {$b_{k}$};
			\node at(4,3) {$b_{2}$};
			\node at(10.5,8) {$a_{k-1}$};
			\node at(4.5,7) {\Large $F$};
			\node at(11,4) {\Large $E$};
		\end{tikzpicture}
	\end{center}
	 \vspace{-1cm} 
	\caption{The boards $E$ and $F$. }\label{fig:board1}
\end{figure}
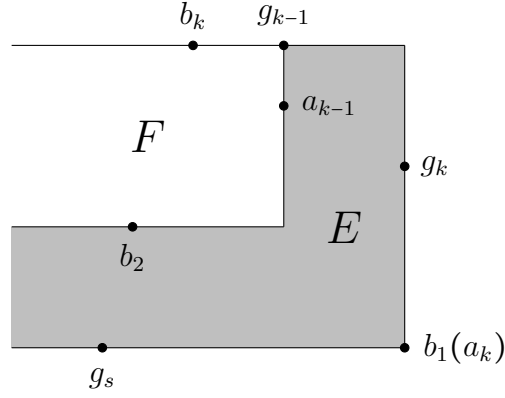

Let $E$ denote the board consisting of (1) the squares that are  below $b_k$ but not below  $b_2$,  and to the right of  $a_{k-1}$ but to the left of $b_1$, and (2) the squares  that above $b_1$ but below $b_2$ and to the left of $b_1$; see Figure \ref{fig:board1}  for an illustration.   
Let $F$ denote the board consisting of the  squares that are   weakly below  $b_{k}$ but weakly above   $b_2$,  and weakly to the left of $a_{k-1}$ as shown in Figure  \ref{fig:board1}. 

For two submatrices $M$  and $M'$ in a transversal $T$,  we say that $M$ is {\em  higher} than $M'$ and $M'$ is {\em lower} than $M$ in $T$ if 
\begin{itemize}
	\item either the lowest square  of $M$ is above the lowest square of $M'$;
	\item or the lowest square of $M$ is located at the same row as  that of $M'$, and the highest square of $M$ is above  that of $M'$. 
\end{itemize}

\begin{lemma}\label{Fact0}
	   There are exactly $k-2$ $1's$ of $\phi^{t}(T)$ located in $F$ and there are no $1's$ of  $\phi^{t}(T)$ located inside $E$.
\end{lemma}
\pf According to the selections of $b_2, b_3, \ldots, b_k$, there are exactly $k-1$ $1's$ of $\phi^{t-1}(T)$  located in  $F$ and there are no $1's$ of  $\phi^{t-1}(T)$  located inside $E$.  By the construction of  $\theta(G)$,  it is easily seen that   $g_s$ is located at the same row as   $b_1$ for some $1\leq s\leq k-2$ in $\phi^t(T)$,   $g_{k-1}$ is located at the same row as $b_k$ and at the same column as $a_{k-1}$, and $g_k$ is positioned at the same column as $b_1$; see    Figure  \ref{fig:board1} for an illustration.
 Since the transformation $\phi$ preserves the positions of the $1's$ of $\phi^{t-1}(T)$ not contained in $G$,  there are exactly $k-2$ $1's$ of $\phi^{t}(T)$ located in $F$ and there are no $1's$ of  $\phi^{t}(T)$ located inside $E$. This completes the proof. \qed

\begin{lemma}\label{lemphi1}
If $\phi^{t}(T)$ contains  a $Q_k$, then this $Q_k$ must be    higher  than $\theta(G)$.  
\end{lemma}
\pf  If not, let  $H$  be  a  $Q_k$ in $\phi^{t}(T)$ that is not higher than $\theta(G)$.  Label the $1's$ of $H$ by $h_1, h_2, \ldots, h_k$ from left to right.   Since $H$ is not higher than $\theta(G)$ and the lowest square of $\theta(G)$ is located at the same  row as $b_1$,   $h_k$ is located not above  $b_1$. 
  We claim that  $h_k$ is located not below  $b_1$.  Otherwise,   since the transformation $\phi$ does not change the positions of other $1's$, at least one of $g_1, g_2, \ldots, g_k$ must fall in $H$ according to the selection of $b_1$.  
 Then  we can get a $Q_k$ in $T$ by replacing $g_i$ with $a_i$ if $g_i$ falls in $H$ for all $1\leq i\leq k$. This yields a contradiction with the selection of $b_1$. Hence, we conclude that  $h_k$  is located at the same row as $b_1$.

   By the construction of  $\theta(G)$,  it is easily seen that   $g_s$ is located at the same row as   $b_1$ for some $1\leq s\leq k-2$ in $\phi^t(T)$,   $g_{k-1}$ is located at the same row as $b_k$ and at the same column as $a_{k-1}$, and $g_k$ is positioned at the same column as $b_1$ as shown in Figure  \ref{fig:board1}.
  This implies that  $h_k=g_s$.  Lemma \ref{Fact0} tells us  that there are exactly $k-2$ $1's$ of $\phi^{t}(T)$ located in $F$  and there are no $1's$ of  $\phi^{t}(T)$ located inside $E$.    
   Therefore,  the topmost $1$ of $H$ must be above $b_k$, which yields a contradiction with the assumption that $H$ is not higher than $\theta(G)$. This completes the proof.  \qed

\begin{lemma}\label{lemphi2}
	   If  $\phi^{t-1}(T)$ contains no  $P_k$  that is higher than $G$, then $\phi^{t}(T)$ contains no $P_k$ that is higher than $\theta(G)$. 
	\end{lemma}
 
	\pf If not, suppose that $H$ is a such $P_k$ in $\phi^t(T)$. Label the $1's$ of $H$ by $h_1, h_2, \ldots, h_k$ from left to right.  Recall that     $g_s$ is located at the same row as   $b_1$ for some $1\leq s\leq k-2$,  $g_{k-1}$ is located at the same row as $b_k$ and at the same column as $a_{k-1}$, and $g_k$ is located at the same column as $b_1$ as shown in Figure  \ref{fig:board1}.  Since   the transformation $\phi$ does not change the positions of other $1's$, at least one of $g_1, g_2, \ldots, g_k$ must fall in $H$.  We shall replace some $1's$ of $H$  to form a $P_k$ in $\phi^{t-1}(T)$ which will contradict the hypothesis.
	 We have two cases.

	\noindent{\bf Case 1.}   The lowest square of $H$ is above  $g_s$ ($b_1$).

 \noindent{\bf Subcase 1.1.}    $h_{k-1}$ is located between $g_k$ and $g_{k-1}$. Lemma \ref{Fact0} states that there are no $1's$ of $\phi^t(T)$ inside $E$. This implies that $h_{k-1}$ must be above $b_k$. In this case,     the $1's$ located in the squares  $b_1, b_2, b_3, \ldots, b_{k-1}$ together with  $h_{k-1}$ would form a $P_k$ in $\phi^{t-1}(T)$ which is higher than $G$,  a contradiction.

  \noindent{\bf Subcase 1.2.}    $h_{k-1}$ is located  weakly to the left of  $g_{k-1}$. Lemma \ref{Fact0} tells us that there are exactly $k-2$ $1's$ of $\phi^{t}(T)$ located in $F$ and there are no $1's$ of  $\phi^{t}(T)$ located inside $E$.  This implies that $h_{k-1}$  must be above $b_k$.   
 Then  we can get a  $P_k$ in $T$ by replacing each $g_i$ with $a_i$ if and only if $g_i$ falls in $H$ for all $1\leq i\leq k-2$ and replacing $h_k$ with $a_k$. It is obvious that such a $P_k$ is higher than $G$ in $\phi^{t-1}(T)$, a contradiction.

    \noindent{\bf Subcase 1.3.}    $h_{k-1}$ is located weakly to the right of  $g_{k}$. Recall that there are exactly $k-2$ $1's$ of $\phi^{t}(T)$ located in $F$ and  there are no $1's$ of  $\phi^{t}(T)$ located inside $E$ by Lemma \ref{Fact0}. Since $g_{k-1}$ is above $g_k$, there are at most $k-3$  $1's$ of $\phi^{t}(T)$ that are  located in $F$ and  are below and to the left of $g_k$.   This implies that $h_{k-1}\neq g_{k}$.        
    By the construction of $\theta(G)$,  for each $g_i$ in $\phi(T)$,  there exists  a $1$ at the same row as $g_i$ in $\phi^{t-1}(T)$ for all $1\leq i\leq k$.  Assume that $a_{j_i}$ is located at the same row   as $g_i$. 
        Then we can  get a  $P_k$ in $T$ by replacing each $g_i$ with   $a_{j_i}$     
                 if and only if $g_i$ falls in $H$. Then such a $P_k$  is higher than $G$ in $\phi^{t-1}(T)$ since $H$ is higher  than $\theta(G)$,   a contradiction.

{ \noindent\bf Case 2. } The lowest $1$ of $H$ is $g_s$.  In this case, $h_{k-1}$ must be  located above $b_k$ ($g_{k-1}$) since $H$ is higher than $\theta(G)$.   

\noindent{\bf Subcase 2.1.}    $h_{k-1}$ is located  weakly to the left of $g_k$.   In this case, 
  we can get a $P_k$ by replacing each $g_i$ with $a_i$ if $g_i$ falls in $H$.  It is apparent that such   a $P_k$ is higher than $G$   in  $\phi^{t-1}(T)$ since  the lowest square of this  $P_k$ is not below $b_1$ and  the topmost $1$ of  this $P_k$  is above $b_k$,   a contradiction. 

\noindent{\bf Subcase 2.2.}  $h_{k-1}$ is located to the right of $g_{k}$. By similar arguments as in Subcase 1.3,    we  can  get a  $P_k$ in $\phi^{t-1}(T)$ by replacing each $g_i$ with   $a_{j_i}$     
if and only if $g_i$ falls in $H$, where $a_{j_i}$ and $g_i$ are located at the same row. Such a $P_k$ is higher than $G$ in $T$  since the lowest $1$ of this $P_k$ is located at the same row as   $b_1$ and the topmost $1$ of  this $P_k$ is located  above $b_k$,   a contradition.  This completes the proof. \qed

\begin{lemma}\label{lemphi3}
 Let $s\geq 1$. Suppose that at the $s$-th application of $\phi$ to $T$  the selected submatrix isomorphic to $Q_k$ is given by $H$. Then  $\phi^{s}(T)$ contains no $P_k$ that is higher than $\theta(H)$. 
\end{lemma}
\pf We prove the assertion by induction on $s$. Suppose that when we apply the transformation $\phi$ to $T$, the selected  submatrix isomorphic to $Q_k$ is given by $H'$. Since $T$ avoids the pattern $P_k$,  $T$ contains no $P_k$  that is higher than $H'$. By Lemma \ref{lemphi2},   $\phi(T)$ contains no $P_k$ that is higher than $\theta(H')$. Hence, the assertion holds for $s=1$.

 Suppose that at the $(s-1)$-th application of  the transformation $\phi$ to $T$, the selected  submatrix isomorphic to $Q_k$ is given by $H''$.  Assume that the assertion holds for $s-1$, that is,  $\phi^{s-1}(T)$ contains no $P_k$ that is higher than $\theta(H'')$.   By Lemma \ref{lemphi1},  the submatrix $H$ isomorphic to $Q_k$ that we select   when we apply   $\phi$ to $\phi^{s-1}(T)$  must be  higher than    $\theta(H'')$ in $\phi^{s-1}(T)$. By the  induction hypothesis,  $\phi^{s-1}(T)$ does not contain a $P_k$ that is higher than $\theta(H'')$. Thus, $\phi^{s-1}(T)$ does not contain a $P_k$ that is higher than $H$.  Then by Lemma  \ref{lemphi2}, we deduce that $\phi^{s}(T)$  does not contain a $P_k$ that is higher than $\theta(H)$  as desired, completing the proof.  \qed

 \begin{lemma}\label{lemphi4}
	Let $s\geq 1$. Suppose that at the $s$-th application of $\phi$ to $T$  the selected submatrix isomorphic to $Q_k$ is given by $H$. Then  $H$ is not isomorphic to $P_k$. 
\end{lemma}
\pf   Clearly, the  assertion holds for $s=1$. For $s\geq 2$, assume that at the $(s-1)$-th application of $\phi$ to $T$, the selected submatrix isomorphic to $Q_k$ is given by $H'$. By Lemma \ref{lemphi1}, the selected submatrix $H$ at the   application of $\phi$ to $\phi^{s-1}(T)$ must be higher than $\theta(H')$.  By Lemma \ref{lemphi3}, $\phi^{s-1}(T)$ contains no $P_k$ that is higher than $\theta(H')$.  This implies that $H$ is not isomorphic to $P_k$ as desired, completing the proof.  \qed

Assume that at the $(s-1)$-th and the   $s$-th  applications of $\phi$ to a transversal $T\in \mathcal{S}_{\lambda}(P_k)$, the selected $Q_k's$  are given by $G$ and $G'$, respectively. 
Lemma  \ref{lemphi1} tells that   
\begin{itemize}
	\item either the lowest square of $G'$ is above  that of $G$,
	\item or the lowest  square of $G'$ is located at the same row as  that of $G$ and the highest  square of $G'$ is above that of $G$.
\end{itemize}
This implies that either the lowest square   or the highest square of the selected $Q_k$ will only go up. Therefore,   our algorithm  $\Phi$ would terminate. Equivalently,  
after applying finitely  many iterations of $\phi$ to a  transversal  $T$ in $\mathcal{S}_{\lambda}(P_k)$, we will eventually  get a transversal $\Phi(T)$ in $\mathcal{S}_{\lambda}(Q_k)$ and hence the map $\Phi$ is well-defined.

	\subsection{The map $\Psi$}
	This subsection is devoted to the construction of the map $\Psi$. To this end, we first    define a transformation $\theta$  on the submatrices of the transversals which changes an occurrence  of $P_k$ to  an occurrence  of $Q_k$.

 \noindent{\bf The transformation $\theta'$.}\\
 Given a transversal $T\in \mathcal{S}_{\lambda}$,  
 let $G$ be a  submatrix  in $T$ that is  isomorphic to  $P_k$ at rows $r_1<r_2<\cdots <r_k$, and columns $c_1<c_2<\cdots<c_k$. Label the $1's$ of $G$ from left to right by $a_1, a_2, \ldots, a_k$, that is,  the $1$ positioned at the column $c_i$ is labeled by $a_i$ for all $1\leq i\leq k$.      Now we generate a  submatrix $\theta'(G)$ isomorphic to $Q_k$  with the same squares of $G$  by   distinguishing the following two cases. 
 \begin{itemize}
 	\item[Case \upshape(\rmnum{1}). ] If $a_{k}$ is  located at row $r_{2}$, 	
 	   we cyclically shift the $1$ located at row $r_i$ to  row $r_{i-2}$ for all $1\leq i\leq k$ with the convention that $r_{-1}:=r_{k-1}$ and $r_{0}:=r_k$. Let $\theta'(G)$ denote the resulting  submatrix. See Figure \ref{fig:theta'1} for an illustration. 
 	\item[Case \upshape(\rmnum{2}).] Otherwise, we cyclically shift the $1$ located at  row $r_i$ to row  $r_{i-1}$ for all $1\leq i\leq k$ with the convention that $r_{0}:=r_k$.  Denote by $G'$  the resulting  submatrix. Let $\theta'(G)$ denote the  submatrix isomorphic to $Q_k$ with the same squares of $G'$ obtained from $G'$ by swapping  the rightmost two columns  and leaving all the other squares fixed.  See Figure \ref{fig:theta'2} for  an example. 
 \end{itemize}
 
 \begin{figure} [H]
 \begin{center}
 	\begin{tikzpicture}[scale = 0.4]
 		\draw (0,0)rectangle(1,1);
 		\draw (1,0)rectangle(2,1);
 		\draw (2,0)rectangle(3,1);
 		\draw (3,0)rectangle(4,1);
 		\draw (0,1)rectangle(1,2);
 		\draw (1,1)rectangle(2,2);
 		\draw (2,1)rectangle(3,2);
 		\draw (3,1)rectangle(4,2);
 		\draw (0,2)rectangle(1,3);
 		\draw (1,2)rectangle(2,3);
 		\draw (2,2)rectangle(3,3);
 		\draw (3,2)rectangle(4,3);
 		\draw (4,2)rectangle(5,3);
 		\draw (5,2)rectangle(6,3);
 		\draw (0,3)rectangle(1,4);
 		\draw (1,3)rectangle(2,4);
 		\draw (2,3)rectangle(3,4);
 		\draw (3,3)rectangle(4,4);
 		\draw (4,3)rectangle(5,4);
 		\draw (5,3)rectangle(6,4);
 		\draw (0,4)rectangle(1,5);
 		\draw (1,4)rectangle(2,5);
 		\draw (2,4)rectangle(3,5);
 		\draw (3,4)rectangle(4,5);
 		\draw (4,4)rectangle(5,5);
 		\draw (5,4)rectangle(6,5);
 		\draw (0,5)rectangle(1,6);
 		\draw (1,5)rectangle(2,6);
 		\draw (2,5)rectangle(3,6);
 		\draw (3,5)rectangle(4,6);
 		\draw (4,5)rectangle(5,6);
 		\draw (5,5)rectangle(6,6);
 		\draw (6,5)rectangle(7,6);
 		\draw (7,5)rectangle(8,6);
 		\draw (0,6)rectangle(1,7);
 		\draw (1,6)rectangle(2,7);
 		\draw (2,6)rectangle(3,7);
 		\draw (3,6)rectangle(4,7);
 		\draw (4,6)rectangle(5,7);
 		\draw (5,6)rectangle(6,7);
 		\draw (6,6)rectangle(7,7);
 		\draw (7,6)rectangle(8,7);
 		\draw (0,7)rectangle(1,8);
 		\draw (1,7)rectangle(2,8);
 		\draw (2,7)rectangle(3,8);
 		\draw (3,7)rectangle(4,8);
 		\draw (4,7)rectangle(5,8);
 		\draw (5,7)rectangle(6,8);
 		\draw (6,7)rectangle(7,8);
 		\draw (7,7)rectangle(8,8);
 		\filldraw[red](0.5,2.5)circle(3pt);
 		\filldraw[black](1.5,0.5)circle(3pt);
 		\filldraw[red](2.5,3.5)circle(3pt);
 		\filldraw[black](3.5,1.5)circle(3pt);
 		\filldraw[red](4.5,7.5)circle(3pt);
 		\filldraw[red](5.5,4.5)circle(3pt);
 		\filldraw[black](6.5,6.5)circle(3pt);
 		\filldraw[black](7.5,5.5)circle(3pt);
 		\node[red] at(0.5,8.5){\footnotesize $c_{1}$};
 		\node[red] at(2.5,8.5){\footnotesize $c_{2}$};
 		\node[red] at(4.5,8.5){\footnotesize $c_{3}$};
 		\node[red] at(5.5,8.5){\footnotesize $c_{4}$};
 		\node[red] at(8.5,7.5){\footnotesize $r_{1}$};
 		\node[red] at(6.5,4.5){\footnotesize $r_{2}$};
 		\node[red] at(6.5,3.5){\footnotesize $r_{3}$};
 		\node[red] at(6.5,2.5){\footnotesize $r_{4}$};
 		
 		\draw [->](9.3,4)--(11.5,4);
 		
 		\draw (13,0)rectangle(14,1);
 		\draw (14,0)rectangle(15,1);
 		\draw (15,0)rectangle(16,1);
 		\draw (16,0)rectangle(17,1);
 		\draw (13,1)rectangle(14,2);
 		\draw (14,1)rectangle(15,2);
 		\draw (15,1)rectangle(16,2);
 		\draw (16,1)rectangle(17,2);
 		\draw (13,2)rectangle(14,3);
 		\draw (14,2)rectangle(15,3);
 		\draw (15,2)rectangle(16,3);
 		\draw (16,2)rectangle(17,3);
 		\draw (17,2)rectangle(18,3);
 		\draw (18,2)rectangle(19,3);
 		\draw (13,3)rectangle(14,4);
 		\draw (14,3)rectangle(15,4);
 		\draw (15,3)rectangle(16,4);
 		\draw (16,3)rectangle(17,4);
 		\draw (17,3)rectangle(18,4);
 		\draw (18,3)rectangle(19,4);
 		\draw (13,4)rectangle(14,5);
 		\draw (14,4)rectangle(15,5);
 		\draw (15,4)rectangle(16,5);
 		\draw (16,4)rectangle(17,5);
 		\draw (17,4)rectangle(18,5);
 		\draw (18,4)rectangle(19,5);
 		\draw (13,5)rectangle(14,6);
 		\draw (14,5)rectangle(15,6);
 		\draw (15,5)rectangle(16,6);
 		\draw (16,5)rectangle(17,6);
 		\draw (17,5)rectangle(18,6);
 		\draw (18,5)rectangle(19,6);
 		\draw (19,5)rectangle(20,6);
 		\draw (20,5)rectangle(21,6);
 		\draw (13,6)rectangle(14,7);
 		\draw (14,6)rectangle(15,7);
 		\draw (15,6)rectangle(16,7);
 		\draw (16,6)rectangle(17,7);
 		\draw (17,6)rectangle(18,7);
 		\draw (18,6)rectangle(19,7);
 		\draw (19,6)rectangle(20,7);
 		\draw (20,6)rectangle(21,7);
 		\draw (13,7)rectangle(14,8);
 		\draw (14,7)rectangle(15,8);
 		\draw (15,7)rectangle(16,8);
 		\draw (16,7)rectangle(17,8);
 		\draw (17,7)rectangle(18,8);
 		\draw (18,7)rectangle(19,8);
 		\draw (19,7)rectangle(20,8);
 		\draw (20,7)rectangle(21,8);
 		\filldraw[red](13.5,4.5)circle(3pt);
 		\filldraw[black](14.5,0.5)circle(3pt);
 		\filldraw[red](15.5,7.5)circle(3pt);
 		\filldraw[black](16.5,1.5)circle(3pt);
 		\filldraw[red](17.5,3.5)circle(3pt);
 		\filldraw[red](18.5,2.5)circle(3pt);
 		\filldraw[black](19.5,6.5)circle(3pt);
 		\filldraw[black](20.5,5.5)circle(3pt);
 		\node[red] at(13.5,8.5){\footnotesize $c_{1}$};
 		\node[red] at(15.5,8.5){\footnotesize $c_{2}$};
 		\node[red] at(17.5,8.5){\footnotesize $c_{3}$};
 		\node[red] at(18.5,8.5){\footnotesize $c_{4}$};
 		\node[red] at(21.5,7.5){\footnotesize $r_{1}$};
 		\node[red] at(19.5,4.5){\footnotesize $r_{2}$};
 		\node[red] at(19.5,3.5){\footnotesize $r_{3}$};
 		\node[red] at(19.5,2.5){\footnotesize $r_{4}$};
 		\node[black] at(4,-2){\footnotesize $G$};
 		\node[black] at(17,-2){\footnotesize $\theta'(G)$};
 		
 	\end{tikzpicture}
 \end{center}
  \vspace{-1cm} 
  \caption{An example of Case \upshape(\rmnum{1}) of the transformation $\theta'$. } \label{fig:theta'1}
 \end{figure}
 
 \begin{figure} [H]
 \begin{center}
 	\begin{tikzpicture}[scale = 0.4]
 		\draw (0,0)rectangle(1,1);
 		\draw (1,0)rectangle(2,1);
 		\draw (2,0)rectangle(3,1);
 		\draw (3,0)rectangle(4,1);
 		\draw (0,1)rectangle(1,2);
 		\draw (1,1)rectangle(2,2);
 		\draw (2,1)rectangle(3,2);
 		\draw (3,1)rectangle(4,2);
 		\draw (0,2)rectangle(1,3);
 		\draw (1,2)rectangle(2,3);
 		\draw (2,2)rectangle(3,3);
 		\draw (3,2)rectangle(4,3);
 		\draw (4,2)rectangle(5,3);
 		\draw (5,2)rectangle(6,3);
 		\draw (0,3)rectangle(1,4);
 		\draw (1,3)rectangle(2,4);
 		\draw (2,3)rectangle(3,4);
 		\draw (3,3)rectangle(4,4);
 		\draw (4,3)rectangle(5,4);
 		\draw (5,3)rectangle(6,4);
 		\draw (0,4)rectangle(1,5);
 		\draw (1,4)rectangle(2,5);
 		\draw (2,4)rectangle(3,5);
 		\draw (3,4)rectangle(4,5);
 		\draw (4,4)rectangle(5,5);
 		\draw (5,4)rectangle(6,5);
 		\draw (0,5)rectangle(1,6);
 		\draw (1,5)rectangle(2,6);
 		\draw (2,5)rectangle(3,6);
 		\draw (3,5)rectangle(4,6);
 		\draw (4,5)rectangle(5,6);
 		\draw (5,5)rectangle(6,6);
 		\draw (6,5)rectangle(7,6);
 		\draw (7,5)rectangle(8,6);
 		\draw (0,6)rectangle(1,7);
 		\draw (1,6)rectangle(2,7);
 		\draw (2,6)rectangle(3,7);
 		\draw (3,6)rectangle(4,7);
 		\draw (4,6)rectangle(5,7);
 		\draw (5,6)rectangle(6,7);
 		\draw (6,6)rectangle(7,7);
 		\draw (7,6)rectangle(8,7);
 		\draw (0,7)rectangle(1,8);
 		\draw (1,7)rectangle(2,8);
 		\draw (2,7)rectangle(3,8);
 		\draw (3,7)rectangle(4,8);
 		\draw (4,7)rectangle(5,8);
 		\draw (5,7)rectangle(6,8);
 		\draw (6,7)rectangle(7,8);
 		\draw (7,7)rectangle(8,8);
 		\filldraw[red](0.5,2.5)circle(3pt);
 		\filldraw[black](1.5,0.5)circle(3pt);
 		\filldraw[red](2.5,4.5)circle(3pt);
 		\filldraw[black](3.5,1.5)circle(3pt);
 		\filldraw[red](4.5,7.5)circle(3pt);
 		\filldraw[red](5.5,3.5)circle(3pt);
 		\filldraw[black](6.5,6.5)circle(3pt);
 		\filldraw[black](7.5,5.5)circle(3pt);
 		\node[red] at(0.5,8.5){\footnotesize $c_{1}$};
 		\node[red] at(2.5,8.5){\footnotesize $c_{2}$};
 		\node[red] at(4.5,8.5){\footnotesize $c_{3}$};
 		\node[red] at(5.5,8.5){\footnotesize $c_{4}$};
 		\node[red] at(8.5,7.5){\footnotesize $r_{1}$};
 		\node[red] at(6.5,4.5){\footnotesize $r_{2}$};
 		\node[red] at(6.5,3.5){\footnotesize $r_{3}$};
 		\node[red] at(6.5,2.5){\footnotesize $r_{4}$};
 		
 		\draw [->](9,4)--(11.5,4);
 		
 		\draw (13,0)rectangle(14,1);
 		\draw (14,0)rectangle(15,1);
 		\draw (15,0)rectangle(16,1);
 		\draw (16,0)rectangle(17,1);
 		\draw (13,1)rectangle(14,2);
 		\draw (14,1)rectangle(15,2);
 		\draw (15,1)rectangle(16,2);
 		\draw (16,1)rectangle(17,2);
 		\draw (13,2)rectangle(14,3);
 		\draw (14,2)rectangle(15,3);
 		\draw (15,2)rectangle(16,3);
 		\draw (16,2)rectangle(17,3);
 		\draw (17,2)rectangle(18,3);
 		\draw (18,2)rectangle(19,3);
 		\draw (13,3)rectangle(14,4);
 		\draw (14,3)rectangle(15,4);
 		\draw (15,3)rectangle(16,4);
 		\draw (16,3)rectangle(17,4);
 		\draw (17,3)rectangle(18,4);
 		\draw (18,3)rectangle(19,4);
 		\draw (13,4)rectangle(14,5);
 		\draw (14,4)rectangle(15,5);
 		\draw (15,4)rectangle(16,5);
 		\draw (16,4)rectangle(17,5);
 		\draw (17,4)rectangle(18,5);
 		\draw (18,4)rectangle(19,5);
 		\draw (13,5)rectangle(14,6);
 		\draw (14,5)rectangle(15,6);
 		\draw (15,5)rectangle(16,6);
 		\draw (16,5)rectangle(17,6);
 		\draw (17,5)rectangle(18,6);
 		\draw (18,5)rectangle(19,6);
 		\draw (19,5)rectangle(20,6);
 		\draw (20,5)rectangle(21,6);
 		\draw (13,6)rectangle(14,7);
 		\draw (14,6)rectangle(15,7);
 		\draw (15,6)rectangle(16,7);
 		\draw (16,6)rectangle(17,7);
 		\draw (17,6)rectangle(18,7);
 		\draw (18,6)rectangle(19,7);
 		\draw (19,6)rectangle(20,7);
 		\draw (20,6)rectangle(21,7);
 		\draw (13,7)rectangle(14,8);
 		\draw (14,7)rectangle(15,8);
 		\draw (15,7)rectangle(16,8);
 		\draw (16,7)rectangle(17,8);
 		\draw (17,7)rectangle(18,8);
 		\draw (18,7)rectangle(19,8);
 		\draw (19,7)rectangle(20,8);
 		\draw (20,7)rectangle(21,8);
 		\filldraw[red](13.5,3.5)circle(3pt);
 		\filldraw[black](14.5,0.5)circle(3pt);
 		\filldraw[red](15.5,7.5)circle(3pt);
 		\filldraw[black](16.5,1.5)circle(3pt);
 		\filldraw[red](17.5,2.5)circle(3pt);
 		\filldraw[red](18.5,4.5)circle(3pt);
 		\filldraw[black](19.5,6.5)circle(3pt);
 		\filldraw[black](20.5,5.5)circle(3pt);
 		\node[red] at(13.5,8.5){\footnotesize $c_{1}$};
 		\node[red] at(15.5,8.5){\footnotesize $c_{2}$};
 		\node[red] at(17.5,8.5){\footnotesize $c_{3}$};
 		\node[red] at(18.5,8.5){\footnotesize $c_{4}$};
 		\node[red] at(21.5,7.5){\footnotesize $r_{1}$};
 		\node[red] at(19.5,4.5){\footnotesize $r_{2}$};
 		\node[red] at(19.5,3.5){\footnotesize $r_{3}$};
 		\node[red] at(19.5,2.5){\footnotesize $r_{4}$};
 		
 		\draw [->](22,4)--(24.5,4);
 		
 		\draw (26,0)rectangle(27,1);
 		\draw (27,0)rectangle(28,1);
 		\draw (28,0)rectangle(29,1);
 		\draw (29,0)rectangle(30,1);
 		\draw (26,1)rectangle(27,2);
 		\draw (27,1)rectangle(28,2);
 		\draw (28,1)rectangle(29,2);
 		\draw (29,1)rectangle(30,2);
 		\draw (26,2)rectangle(27,3);
 		\draw (27,2)rectangle(28,3);
 		\draw (28,2)rectangle(29,3);
 		\draw (29,2)rectangle(30,3);
 		\draw (30,2)rectangle(31,3);
 		\draw (31,2)rectangle(32,3);
 		\draw (26,3)rectangle(27,4);
 		\draw (27,3)rectangle(28,4);
 		\draw (28,3)rectangle(29,4);
 		\draw (29,3)rectangle(30,4);
 		\draw (30,3)rectangle(31,4);
 		\draw (31,3)rectangle(32,4);
 		\draw (26,4)rectangle(27,5);
 		\draw (27,4)rectangle(28,5);
 		\draw (28,4)rectangle(29,5);
 		\draw (29,4)rectangle(30,5);
 		\draw (30,4)rectangle(31,5);
 		\draw (31,4)rectangle(32,5);
 		\draw (26,5)rectangle(27,6);
 		\draw (27,5)rectangle(28,6);
 		\draw (28,5)rectangle(29,6);
 		\draw (29,5)rectangle(30,6);
 		\draw (30,5)rectangle(31,6);
 		\draw (31,5)rectangle(32,6);
 		\draw (32,5)rectangle(33,6);
 		\draw (33,5)rectangle(34,6);
 		\draw (26,6)rectangle(27,7);
 		\draw (27,6)rectangle(28,7);
 		\draw (28,6)rectangle(29,7);
 		\draw (29,6)rectangle(30,7);
 		\draw (30,6)rectangle(31,7);
 		\draw (31,6)rectangle(32,7);
 		\draw (32,6)rectangle(33,7);
 		\draw (33,6)rectangle(34,7);
 		\draw (26,7)rectangle(27,8);
 		\draw (27,7)rectangle(28,8);
 		\draw (28,7)rectangle(29,8);
 		\draw (29,7)rectangle(30,8);
 		\draw (30,7)rectangle(31,8);
 		\draw (31,7)rectangle(32,8);
 		\draw (32,7)rectangle(33,8);
 		\draw (33,7)rectangle(34,8);
 		\filldraw[red](26.5,3.5)circle(3pt);
 		\filldraw[black](27.5,0.5)circle(3pt);
 		\filldraw[red](28.5,7.5)circle(3pt);
 		\filldraw[black](29.5,1.5)circle(3pt);
 		\filldraw[red](30.5,4.5)circle(3pt);
 		\filldraw[red](31.5,2.5)circle(3pt);
 		\filldraw[black](32.5,6.5)circle(3pt);
 		\filldraw[black](33.5,5.5)circle(3pt);
 		\node[red] at(26.5,8.5){\footnotesize $c_{1}$};
 		\node[red] at(28.5,8.5){\footnotesize $c_{2}$};
 		\node[red] at(30.5,8.5){\footnotesize $c_{3}$};
 		\node[red] at(31.5,8.5){\footnotesize $c_{4}$};
 		\node[red] at(34.5,7.5){\footnotesize $r_{1}$};
 		\node[red] at(32.5,4.5){\footnotesize $r_{2}$};
 		\node[red] at(32.5,3.5){\footnotesize $r_{3}$};
 		\node[red] at(32.5,2.5){\footnotesize $r_{4}$};
 		\node[black] at(4,-2){\footnotesize $G$};
 		\node[black] at(17,-2){\footnotesize $G'$};
 		\node[black] at(30,-2){\footnotesize $\theta'(G)$};
 	\end{tikzpicture}
 \end{center}
  \vspace{-0.8cm} 
  \caption{An example of Case \upshape(\rmnum{2}) of the transformation $\theta'$. } \label{fig:theta'2}
 \end{figure}
 
 The following two observations follow directly from the definitions of the transformations $\theta$ and $\theta'$. 
 \begin{observation}\label{ob1}
 	Let $G$ be a submatrix isomorphic to $Q_k$ in which the $1's$ are labeled by $a_1, a_2, \ldots, a_k$ from left to right. If  $a_{k-1}$ is not positioned at the topmost row of $G$, then the procedures in Case (\rmnum{1}) and  Case (\rmnum{2}) of $\theta'$ reverse  the procedures in Case (\rmnum{1}) and  Case (\rmnum{2}) of $\theta$, respectively. Hence, we have $\theta'(\theta (G))=G$ when  $a_{k-1}$ is not positioned at the topmost row of $G$. 
 	\end{observation}

  \begin{observation}\label{ob2}
 	Let $G$ be a submatrix isomorphic to $P_k$ in which the $1's$ are labeled by $a_1, a_2, \ldots, a_k$ from left to right. If  $a_k$ is not positioned at the lowest row of $G$, then the procedures in Case (\rmnum{1}) and  Case (\rmnum{2}) of $\theta$ reverse  the procedures in Case (\rmnum{1}) and  Case (\rmnum{2}) of $\theta'$, respectively. Hence, we have $\theta(\theta' (G))=G$ when $a_k$ is not positioned at the lowest row of $G$. 
 \end{observation}
 {\noindent \bf The map  $\Psi: \mathcal{S}_{\lambda}(Q_k)\rightarrow \mathcal{S}_{\lambda}(P_k)$ .}\\
 Let $T$ be a transversal  in  $ \mathcal{S}_{\lambda}(Q_k)$. 
 \begin{itemize}
 	\item[Step 1.] If  $T$ contains no $P_k$, terminate.
 	 	\item[Step 2.]  Find the highest square $b_1$  containing   $1$
 	 	such that there is a   $P_k$ in which  the lowest $1$ of its leftmost $k-1$ $1's$ is located in $b_1$.
 	 	\item[Step 3.] Find the highest  square $b_{k-1}$  containing   $1$
 	 	such that there is a   $P_k$ in which   the lowest $1$ of its leftmost $k-1$ $1's$ is located in $b_1$ and the topmost $1$ is located in $b_{k-1}$. 	
 	 	 
   	\item[Step 4.]  Find the highest  squares  $b_2, \ldots, b_{k-2}$ (listed  from bottom to top) containing  $1's$   that are   located above $b_1$ but below $b_{k-1}$ and   to the left of  $b_{k-1}$. 
   	\item[Step 5.]
   	 Find the leftmost  square $b_k$ containing   $1$ such that there is a submatrix $G$ isomorphic to $P_k$ in $T$ so that  the squares $b_1, b_2, \ldots, b_k$  contain the $1's$ of $G$.  	 
 	\item[Step 6.]  Replace  the submatrix $G$  with $\theta'(G)$ and leave all the other squares fixed. 	 
 	\item[Step 7.]  Repeat the above procedure  until  there is no $P_k$ left. Then we get a $\Psi(T)\in  \mathcal{S}_{\lambda}(P_k)$.
 \end{itemize}

 For our convenience,    Steps 2--6 of the algorithm $\Psi$  is called  {\em the transformation $\psi$}.
 Let  $t$ be a nonnegative  integer. Denote by $\psi^{t}(T)$ the resulting transversal  obtained  by the $t$-th application of $\psi $ to the $T$. Clearly, we have $\psi^{0}(T)=T$. 
 
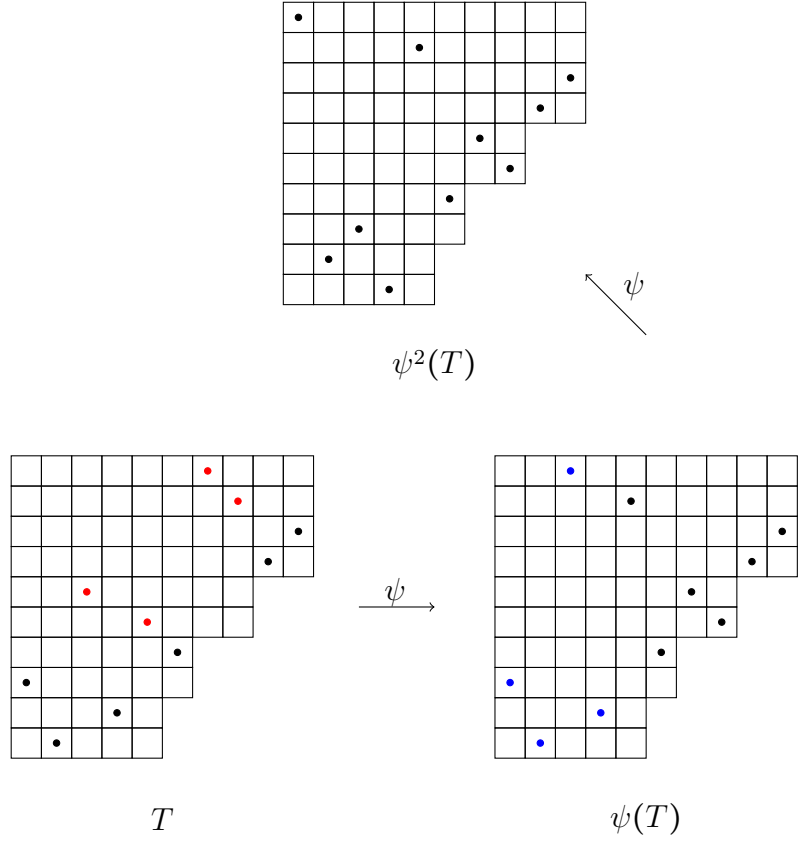
\begin{figure} [H]
 \begin{center}
 	\begin{tikzpicture}[scale = 0.4]
 		\draw (0,0)rectangle(1,1);
 		\draw (1,0)rectangle(2,1);
 		\draw (2,0)rectangle(3,1);
 		\draw (3,0)rectangle(4,1);
 		\draw (4,0)rectangle(5,1);
 		\draw (0,1)rectangle(1,2);
 		\draw (1,1)rectangle(2,2);
 		\draw (2,1)rectangle(3,2);
 		\draw (3,1)rectangle(4,2);
 		\draw (4,1)rectangle(5,2);
 		\draw (0,2)rectangle(1,3);
 		\draw (1,2)rectangle(2,3);
 		\draw (2,2)rectangle(3,3);
 		\draw (3,2)rectangle(4,3);
 		\draw (4,2)rectangle(5,3);
 		\draw (5,2)rectangle(6,3);
 		\draw (0,3)rectangle(1,4);
 		\draw (1,3)rectangle(2,4);
 		\draw (2,3)rectangle(3,4);
 		\draw (3,3)rectangle(4,4);
 		\draw (4,3)rectangle(5,4);
 		\draw (5,3)rectangle(6,4);
 		\draw (0,4)rectangle(1,5);
 		\draw (1,4)rectangle(2,5);
 		\draw (2,4)rectangle(3,5);
 		\draw (3,4)rectangle(4,5);
 		\draw (4,4)rectangle(5,5);
 		\draw (5,4)rectangle(6,5);
 		\draw (6,4)rectangle(7,5);
 		\draw (7,4)rectangle(8,5);
 		\draw (0,5)rectangle(1,6);
 		\draw (1,5)rectangle(2,6);
 		\draw (2,5)rectangle(3,6);
 		\draw (3,5)rectangle(4,6);
 		\draw (4,5)rectangle(5,6);
 		\draw (5,5)rectangle(6,6);
 		\draw (6,5)rectangle(7,6);
 		\draw (7,5)rectangle(8,6);
 		\draw (0,6)rectangle(1,7);
 		\draw (1,6)rectangle(2,7);
 		\draw (2,6)rectangle(3,7);
 		\draw (3,6)rectangle(4,7);
 		\draw (4,6)rectangle(5,7);
 		\draw (5,6)rectangle(6,7);
 		\draw (6,6)rectangle(7,7);
 		\draw (7,6)rectangle(8,7);
 		\draw (8,6)rectangle(9,7);
 		\draw (9,6)rectangle(10,7);
 		\draw (0,7)rectangle(1,8);
 		\draw (1,7)rectangle(2,8);
 		\draw (2,7)rectangle(3,8);
 		\draw (3,7)rectangle(4,8);
 		\draw (4,7)rectangle(5,8);
 		\draw (5,7)rectangle(6,8);
 		\draw (6,7)rectangle(7,8);
 		\draw (7,7)rectangle(8,8);
 		\draw (8,7)rectangle(9,8);
 		\draw (9,7)rectangle(10,8);
 		\draw (0,8)rectangle(1,9);
 		\draw (1,8)rectangle(2,9);
 		\draw (2,8)rectangle(3,9);
 		\draw (3,8)rectangle(4,9);
 		\draw (4,8)rectangle(5,9);
 		\draw (5,8)rectangle(6,9);
 		\draw (6,8)rectangle(7,9);
 		\draw (7,8)rectangle(8,9);
 		\draw (8,8)rectangle(9,9);
 		\draw (9,8)rectangle(10,9);
 		\draw (0,9)rectangle(1,10);
 		\draw (1,9)rectangle(2,10);
 		\draw (2,9)rectangle(3,10);
 		\draw (3,9)rectangle(4,10);
 		\draw (4,9)rectangle(5,10);
 		\draw (5,9)rectangle(6,10);
 		\draw (6,9)rectangle(7,10);
 		\draw (7,9)rectangle(8,10);
 		\draw (8,9)rectangle(9,10);
 		\draw (9,9)rectangle(10,10);
 		\filldraw[black](0.5,2.5)circle(3pt);
 		\filldraw[black](1.5,0.5)circle(3pt);
 		\filldraw[red](2.5,5.5)circle(3pt);
 		\filldraw[black](3.5,1.5)circle(3pt);
 		\filldraw[red](4.5,4.5)circle(3pt);
 		\filldraw[black](5.5,3.5)circle(3pt);
 		\filldraw[red](6.5,9.5)circle(3pt);
 		\filldraw[red](7.5,8.5)circle(3pt);
 		\filldraw[black](8.5,6.5)circle(3pt);
 		\filldraw[black](9.5,7.5)circle(3pt);
 		\node at(5,-2){$T$};
 		
 		\draw [->](11.5,5)--(14,5);
 		\node at(12.7,5.5){$\psi$};
 		
 		\draw (16,0)rectangle(17,1);
 		\draw (17,0)rectangle(18,1);
 		\draw (18,0)rectangle(19,1);
 		\draw (19,0)rectangle(20,1);
 		\draw (20,0)rectangle(21,1);
 		\draw (16,1)rectangle(17,2);
 		\draw (17,1)rectangle(18,2);
 		\draw (18,1)rectangle(19,2);
 		\draw (19,1)rectangle(20,2);
 		\draw (20,1)rectangle(21,2);
 		\draw (16,2)rectangle(17,3);
 		\draw (17,2)rectangle(18,3);
 		\draw (18,2)rectangle(19,3);
 		\draw (19,2)rectangle(20,3);
 		\draw (20,2)rectangle(21,3);
 		\draw (21,2)rectangle(22,3);
 		\draw (16,3)rectangle(17,4);
 		\draw (17,3)rectangle(18,4);
 		\draw (18,3)rectangle(19,4);
 		\draw (19,3)rectangle(20,4);
 		\draw (20,3)rectangle(21,4);
 		\draw (21,3)rectangle(22,4);
 		\draw (16,4)rectangle(17,5);
 		\draw (17,4)rectangle(18,5);
 		\draw (18,4)rectangle(19,5);
 		\draw (19,4)rectangle(20,5);
 		\draw (20,4)rectangle(21,5);
 		\draw (21,4)rectangle(22,5);
 		\draw (22,4)rectangle(23,5);
 		\draw (23,4)rectangle(24,5);
 		\draw (16,5)rectangle(17,6);
 		\draw (17,5)rectangle(18,6);
 		\draw (18,5)rectangle(19,6);
 		\draw (19,5)rectangle(20,6);
 		\draw (20,5)rectangle(21,6);
 		\draw (21,5)rectangle(22,6);
 		\draw (22,5)rectangle(23,6);
 		\draw (23,5)rectangle(24,6);
 		\draw (16,6)rectangle(17,7);
 		\draw (17,6)rectangle(18,7);
 		\draw (18,6)rectangle(19,7);
 		\draw (19,6)rectangle(20,7);
 		\draw (20,6)rectangle(21,7);
 		\draw (21,6)rectangle(22,7);
 		\draw (22,6)rectangle(23,7);
 		\draw (23,6)rectangle(24,7);
 		\draw (24,6)rectangle(25,7);
 		\draw (25,6)rectangle(26,7);
 		\draw (16,7)rectangle(17,8);
 		\draw (17,7)rectangle(18,8);
 		\draw (18,7)rectangle(19,8);
 		\draw (19,7)rectangle(20,8);
 		\draw (20,7)rectangle(21,8);
 		\draw (21,7)rectangle(22,8);
 		\draw (22,7)rectangle(23,8);
 		\draw (23,7)rectangle(24,8);
 		\draw (24,7)rectangle(25,8);
 		\draw (25,7)rectangle(26,8);
 		\draw (16,8)rectangle(17,9);
 		\draw (17,8)rectangle(18,9);
 		\draw (18,8)rectangle(19,9);
 		\draw (19,8)rectangle(20,9);
 		\draw (20,8)rectangle(21,9);
 		\draw (21,8)rectangle(22,9);
 		\draw (22,8)rectangle(23,9);
 		\draw (23,8)rectangle(24,9);
 		\draw (24,8)rectangle(25,9);
 		\draw (25,8)rectangle(26,9);
 		\draw (16,9)rectangle(17,10);
 		\draw (17,9)rectangle(18,10);
 		\draw (18,9)rectangle(19,10);
 		\draw (19,9)rectangle(20,10);
 		\draw (20,9)rectangle(21,10);
 		\draw (21,9)rectangle(22,10);
 		\draw (22,9)rectangle(23,10);
 		\draw (23,9)rectangle(24,10);
 		\draw (24,9)rectangle(25,10);
 		\draw (25,9)rectangle(26,10);
 		\filldraw[blue](16.5,2.5)circle(3pt);
 		\filldraw[blue](17.5,0.5)circle(3pt);
 		\filldraw[blue](18.5,9.5)circle(3pt);
 		\filldraw[blue](19.5,1.5)circle(3pt);
 		\filldraw[black](20.5,8.5)circle(3pt);
 		\filldraw[black](21.5,3.5)circle(3pt);
 		\filldraw[black](22.5,5.5)circle(3pt);
 		\filldraw[black](23.5,4.5)circle(3pt);
 		\filldraw[black](24.5,6.5)circle(3pt);
 		\filldraw[black](25.5,7.5)circle(3pt);
 		\node at(21,-2){$\psi(T)$};
 		
 		\draw [->](21,14)--(19,16);
 		\node at(20.6,15.6){$\psi$};
 		
 		\draw (9,15)rectangle(10,16);
 		\draw (10,15)rectangle(11,16);
 		\draw (11,15)rectangle(12,16);
 		\draw (12,15)rectangle(13,16);
 		\draw (13,15)rectangle(14,16);
 		\draw (9,16)rectangle(10,17);
 		\draw (10,16)rectangle(11,17);
 		\draw (11,16)rectangle(12,17);
 		\draw (12,16)rectangle(13,17);
 		\draw (13,16)rectangle(14,17);
 		\draw (9,17)rectangle(10,18);
 		\draw (10,17)rectangle(11,18);
 		\draw (11,17)rectangle(12,18);
 		\draw (12,17)rectangle(13,18);
 		\draw (13,17)rectangle(14,18);
 		\draw (14,17)rectangle(15,18);
 		\draw (9,18)rectangle(10,19);
 		\draw (10,18)rectangle(11,19);
 		\draw (11,18)rectangle(12,19);
 		\draw (12,18)rectangle(13,19);
 		\draw (13,18)rectangle(14,19);
 		\draw (14,18)rectangle(15,19);
 		\draw (9,19)rectangle(10,20);
 		\draw (10,19)rectangle(11,20);
 		\draw (11,19)rectangle(12,20);
 		\draw (12,19)rectangle(13,20);
 		\draw (13,19)rectangle(14,20);
 		\draw (14,19)rectangle(15,20);
 		\draw (15,19)rectangle(16,20);
 		\draw (16,19)rectangle(17,20);
 		\draw (9,20)rectangle(10,21);
 		\draw (10,20)rectangle(11,21);
 		\draw (11,20)rectangle(12,21);
 		\draw (12,20)rectangle(13,21);
 		\draw (13,20)rectangle(14,21);
 		\draw (14,20)rectangle(15,21);
 		\draw (15,20)rectangle(16,21);
 		\draw (16,20)rectangle(17,21);
 		\draw (9,21)rectangle(10,22);
 		\draw (10,21)rectangle(11,22);
 		\draw (11,21)rectangle(12,22);
 		\draw (12,21)rectangle(13,22);
 		\draw (13,21)rectangle(14,22);
 		\draw (14,21)rectangle(15,22);
 		\draw (15,21)rectangle(16,22);
 		\draw (16,21)rectangle(17,22);
 		\draw (17,21)rectangle(18,22);
 		\draw (18,21)rectangle(19,22);
 		\draw (9,22)rectangle(10,23);
 		\draw (10,22)rectangle(11,23);
 		\draw (11,22)rectangle(12,23);
 		\draw (12,22)rectangle(13,23);
 		\draw (13,22)rectangle(14,23);
 		\draw (14,22)rectangle(15,23);
 		\draw (15,22)rectangle(16,23);
 		\draw (16,22)rectangle(17,23);
 		\draw (17,22)rectangle(18,23);
 		\draw (18,22)rectangle(19,23);
 		\draw (9,23)rectangle(10,24);
 		\draw (10,23)rectangle(11,24);
 		\draw (11,23)rectangle(12,24);
 		\draw (12,23)rectangle(13,24);
 		\draw (13,23)rectangle(14,24);
 		\draw (14,23)rectangle(15,24);
 		\draw (15,23)rectangle(16,24);
 		\draw (16,23)rectangle(17,24);
 		\draw (17,23)rectangle(18,24);
 		\draw (18,23)rectangle(19,24);
 		\draw (9,24)rectangle(10,25);
 		\draw (10,24)rectangle(11,25);
 		\draw (11,24)rectangle(12,25);
 		\draw (12,24)rectangle(13,25);
 		\draw (13,24)rectangle(14,25);
 		\draw (14,24)rectangle(15,25);
 		\draw (15,24)rectangle(16,25);
 		\draw (16,24)rectangle(17,25);
 		\draw (17,24)rectangle(18,25);
 		\draw (18,24)rectangle(19,25);
 		\filldraw[black](9.5,24.5)circle(3pt);
 		\filldraw[black](10.5,16.5)circle(3pt);
 		\filldraw[black](11.5,17.5)circle(3pt);
 		\filldraw[black](12.5,15.5)circle(3pt);
 		\filldraw[black](13.5,23.5)circle(3pt);
 		\filldraw[black](14.5,18.5)circle(3pt);
 		\filldraw[black](15.5,20.5)circle(3pt);
 		\filldraw[black](16.5,19.5)circle(3pt);
 		\filldraw[black](17.5,21.5)circle(3pt);
 		\filldraw[black](18.5,22.5)circle(3pt);
 		\node at(14,13){$\psi^2(T)$};
 		
 	\end{tikzpicture}
  
\end{center}
 \vspace{-0.8cm} 
\caption{An example of the map $\Psi$. }\label{fig:Psi}
\end{figure}

 Let $T$ be the transversal as shown in Figure \ref{fig:Psi}.  
  When we apply the transformation $\psi$ to $T$,  the $1's$ of the selected submatrix $G$ isomorphic to  $P_4$ is colored by red. By replacing $G$ with $\theta'(G)$, we get a transversal $\psi(T)$ displayed in Figure \ref{fig:Psi}. At the second application of $\Psi$ to $T$, the $1's$ of the selected submatrix $G'$ is colored by blue. By replacing $G'$ with $\theta'(G')$, we finally get a transversal $\psi^2(T)$ that avoids $P_4$  as shown in Figure  \ref{fig:Psi}.
  
  Now we proceed to prove that the transformation $\psi$ has the following celebrated properties.  We begin with some necessary     notations.  
  
  If $T$ contains $P_k$,  
  suppose that at the $t$-th application of $\psi$ to $T$, we select the submartx $G$ isomorphic to $P_k$ in which  the squares containing $1's$ are labeled  as the same as those in the definition of $\psi$. For our convenience, we  also
  label     the $1's$ of $G$ by $a_1, a_2, \ldots, a_k$ from left to right and   label the $1's$ of  the  submatrix $\theta(G)$ from left to right by $g_1,g_2, \ldots, g_k$.

\begin{figure} [H]
 	\begin{center}
 		\begin{tikzpicture}[scale = 0.4]
 			
 			\fill[lightgray](29,0)rectangle(33,10);
 			\draw[black](20,0)--(33,0);
 			\draw[black](33,0)--(33,10);
 			\draw[black](20,10)--(33,10);
 			\draw[black](29,10)--(29,0);
 			\filldraw[black](33,0)circle(4pt);
 			\filldraw[black](29,10)circle(4pt);
 			\filldraw[black](24,0)circle(4pt);
 			\filldraw[black](33,4)circle(4pt);
 			\node at(34.5,0) {$x (g_k)$};
 			\node at(29,11) {$b_{k-1}$};
 			\node at(24,-1) {$b_{1}$};
 			\node at(35,4) {$b_{k} (a_k)$};
 			\node at(24.5,5) {\Large $F'$};
 			\node at(31,5) {\Large $E'$};
 				\node at(26,11) {$g_s$};
 				\filldraw[black](26,10)circle(4pt);
 		\end{tikzpicture}
 		
 	\end{center}
 	 \vspace{-0.8cm} 
 	\caption{The boards $E'$ and $F'$. }\label{fig:board}
 \end{figure}
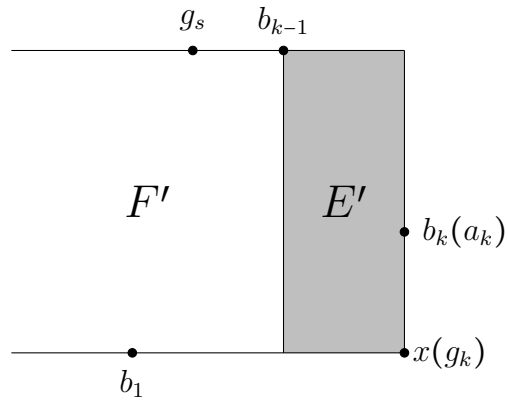

 Let $x$ be  the square located  in  the lowest right corner of the $G$.
 If  $G$ is not isomorphic to $Q_k$, the lowest $1$ of $G$ must be located in $b_1$. Let $E'$ denote the board  consisting of   the squares that are below $b_{k-1}$ but above $b_1$, and to the right of  $b_{k-1}$ but to the left of $b_k$; see Figure  \ref{fig:board}  for an illustration.
 Let $F'$ denote the board consisting of the squares  that are  above $b_1$ but weakly   below $b_{k-1}$,  and weakly to the left of $b_{k-1}$ as shown in Figure  \ref{fig:board}.

 \begin{lemma}\label{Fact1}
 If $G$ is not isomorphic to $Q_k$, then there are no $1's$ of    $\psi^{t}(T)$ inside  $E'$  and there are exactly $k-1$ $1's$ of    $\psi^{t}(T)$ located in  $F'$. 
 \end{lemma}
 \pf  By the selection of $b_1$ and the selection of $b_k$, one can easily check that there are exactly $k-2$ $1's$  positioned  in  $F'$ and  there does not exist any $1$ inside   $E'$  in $\psi^{t-1}(T)$.    By the construction of $\theta'(G)$,  the transformation $\psi$ preserves the positions of $1's$ that are not contained in $G$.  Moreover, $g_k$ is positioned at $x$ and $g_s$ is located  at the same row as $b_{k-1}$ for some $1\leq s\leq k-2$  as illustrated in Figure \ref{fig:board}. 
 Therefore,  there are no $1's$ of    $\psi^{t}(T)$ inside   $E'$ and there are exactly    $k-1$ $1's$  of    $\psi^{t}(T)$  located in    $F'$, completing the proof. \qed

  \begin{lemma}\label{lempsi1}
 	If $G$ is not isomorphic to $Q_k$ and  $\psi^{t-1}(T)$ contains no  $Q_k$  that is lower than $G$, then $\psi^{t}(T)$ contains no $Q_k$ that is lower than $\theta'(G)$. 
 \end{lemma}
 \pf  
      If not, suppose that $H$ is a such $Q_k$ in $\psi^t(T)$. Label the $1's$ of $H$ by $h_1, h_2, \ldots, h_k$ from left to right.  Lemma \ref{Fact1} states that there are no $1's$  of  $\psi^t(T)$ inside $E'$ and there are exactly $k-1$ $1's$ of  $\psi^t(T)$ positioned in   $F'$.          
       This implies that $h_k$ must be below $x$ that is the lowest right corner of $\theta'(G)$.   
  Since   the transformation $\psi$ does not change the positions of other $1's$, at least one of $g_1, g_2, \ldots, g_k$ must fall in $H$. We can obtain a $Q_k$ in $\psi^{t-1}(T)$ from $H$ by replacing each $g_i$ with $a_i$ if $g_i$ falls in $H$. Clearly, such a $Q_k$ is lower than $G$ in  $\psi^{t-1}(T)$, which contradicts the hypothesis. This completes the proof. \qed

  \begin{lemma}\label{lempsi2}
 	If $G$ is not isomorphic to $Q_k$ and $\psi^{t}(T)$ still contains a $P_k$, then such a $P_k$ must be lower   than  $\theta'(G)$. 
 \end{lemma}
 \pf   
  Suppose that $G'$ is   a $P_k$ that is not lower than $\theta'(G)$ in $\psi^t(T)$. Label the  $1's$  of $G'$   by $g'_1, g'_2, \ldots, g'_k$ from left to right.    
  We shall replace some $1's$ of $G'$  to form a $P_k$ in $\psi^{t-1}(T)$ which will either  contradict the selection of $b_1$ or the selection of $b_{k-1}$.    
 We have two cases.

 \noindent{\bf Case 1.} $g'_{k-1}$ is located to the left of $x$.\\
  By the construction of $\theta'(G)$,    $g_k$ is positioned at $x$ and $g_s$ is located  at the same row as $b_{k-1}$ for some $1\leq s\leq k-2$  as illustrated in Figure \ref{fig:board}. 
  Lemma \ref{Fact1} tells us that  there are exactly $k-1$ $1's$ of   $\psi^{t}(T)$ that are positioned  in $F'$, namely, $g_1, g_2, \ldots, g_{k-1}$, and there are no $1's$ of    $\psi^{t}(T)$ that are located inside $E'$. Since $s\leq k-2$,   it follows that $g'_{k-1}$ must be above    $b_{k-1}$  and the lowest $1$ of $G'$ must not be below  $x$. Then we can get a submatrix $G''$ isomorphic to  $P_k$ in $\psi^{t-1}(T)$ from $G'$  by replacing $g_i$ with $a_i$ if $g_i$ falls in $G'$ for all $1\leq i<k$ and replacing $g'_k$ with $a_k$.  Clearly, 
 \begin{itemize}
 	\item either the  lowest $1$ of the first $k-1$ $1's$ of $G''$ is positioned above $b_1$,
 	\item or  the  lowest $1$ of the first $k-1$ $1's$ of $G''$ is positioned at  $b_1$ and the the highest $1$ of $G''$ is above $b_{k-1}$.
 \end{itemize}
 The former would contradict the selection of $b_1$, whereas the latter would contradict the selection of $b_{k-1}$.  
 
 \noindent{\bf Case 2.} $g'_{k-1}$ is located to the right of $x$. \\
 By the construction of $\theta'(G)$,  for each $g_i$ in $\psi^{t}(T)$,  there exists  a $1$ at the same row as $g_i$ in $\psi^{t-1}(T)$ for all $1\leq i\leq k$.  Assume that $a_{j_i}$ is located at the same row   as $g_i$. Then we can get a  $G''$ isomorphic to $P_k$ in $\psi^{t-1}(T)$ by replacing each $g_i$ with $a_{j_i}$ if $g_i$ falls in $G'$.  
 Since $G'$ is not lower $\theta'(G)$ in $\psi^{t}(T)$,   it follows that $G''$ is not lower than $G$ in  $\psi^{t-1}(T)$.  Observe that the highest $1$ of $G''$  is  $g'_{k-1}$.

 \noindent{\bf Subcase 2.1.} 
 $g'_{k-1}$ is above $b_{k-1}$. Since $G''$ is not lower than $G$, the lowest $1$ of $G''$ must be located not below $b_1$. According to the  selection of $b_1$,  one can see that the lowest $1$ of $G''$ must be located at  $b_1$.   Recall that  the highest $1$ of $G''$ is given by $g'_{k-1}$ that is above $b_{k-1}$.   This   yields a contradiction with  the selection of $b_{k-1}$.

 \noindent{\bf Subcase 2.2.} 
 $g'_{k-1}$ is below $b_{k-1}$.  Since $G''$ is not lower than $G$,  the lowest $1$ of $G''$ must be above $b_1$. This  also contradicts the selection of $b_1$,  completing the proof. \qed

 \begin{lemma}\label{lempsi3}
 	Let $s\geq 1$. Suppose that at the $s$-th application of $\psi$ to $T$  the selected submatrix isomorphic to $P_k$ is given by $H$. Then  $\psi^{s}(T)$ contains neither $Q_k$ that is lower than $\theta'(H)$ nor $P_k$ that is not lower than $\theta'(H)$. 
 \end{lemma}
 \pf We prove the assertion by induction on $s$. Suppose that when we apply the transformation $\psi$ to $T$, the selected  submatrix isomorphic to $P_k$ is given by $G$. Since $T$ avoids $Q_k$,  $T$ contains no $Q_k$ that is below $G$ and $G$ is not isomorphic to $Q_k$. Then by Lemmas \ref{lempsi1} and \ref{lempsi2}, it follows that $\psi(T)$ contains  neither    $Q_k$ that is lower than $\theta'(G)$  nor  $P_k$ that is not lower than $\theta'(G)$.
   Therefore,  the assertion holds for $s=1$.

 Suppose that at the $(s-1)$-th application of  the transformation $\psi$ to $T$, the selected  submatrix isomorphic to $P_k$ is given by $H'$.  Assume that the assertion holds for $s-1$, that is,  $\psi^{s-1}(T)$ contains neither $Q_k$ that is lower than $\theta'(H')$ nor $P_k$ that is not lower than $\theta'(H')$.   This implies that the submatrix $H$ isomorphic to $P_k$ that we select   when we apply   $\psi$ to $\psi^{s-1}(T)$  must be  lower than    $\theta'(H')$ in $\psi^{s-1}(T)$. By the  induction hypothesis,  $\psi^{s-1}(T)$   contains no   $Q_k$ that is lower than $\theta'(H')$.  Thus, $\psi^{s-1}(T)$ does not contain a $Q_k$ that is lower than $H$  and $H$ is not isomorphic to $Q_k$. By Lemmas  \ref{lempsi1} and \ref{lempsi2}, we deduce that $\psi^{s}(T)$   contains neither $Q_k$ that is lower than $\theta'(H)$ nor $P_k$ that is not lower than $\theta'(H)$ as desired. This completes the proof.   \qed

  \begin{lemma}\label{lempsi4}
 	Let $s\geq 1$. Suppose that at the $s$-th application of $\psi$ to $T$  the selected submatrix isomorphic to $P_k$ is given by $H$.  Then  $H$ is not isomorphic to $Q_k$. 
 \end{lemma}
 \pf   It is apparent that  the assertion holds for $s=1$. 
 For $s\geq 2$, assume that at the $(s-1)$-th application of $\psi$ to $T$, the selected submatrix isomorphic to $P_k$ is given by $H'$. By Lemma \ref{lempsi3}, the selected submatrix $H$ at the   application of $\psi$ to $\psi^{s-1}(T)$ must be lower than $\theta'(H')$. Again by Lemma \ref{lempsi3}, $\psi^{s-1}(T)$ contains no $Q_k$ that is lower than $\theta'(H')$.  This implies that $H$ is not isomorphic to $Q_k$  as desired, completing the proof.  \qed

 Assume that at the $(s-1)$-th and the   $s$-th  applications of $\psi$ to a transversal $T\in \mathcal{S}_{\lambda}(Q_k)$, the selected $P_k's$  are given by $G$ and $G'$, respectively. 
  Lemma  \ref{lempsi3} tells that   
  \begin{itemize}
  	\item either the lowest square of $G'$ is below  that of $G$,
  	\item or the lowest  square of $G'$ is located at the same row as  that of $G$ and the highest  square of $G'$ is below that of $G$.
  	\end{itemize}
  This implies that  either the lowest square   or the highest square of the selected $P_k$ will  go down.  Therefore,  our algorithm would terminate, that is, 
  after applying finitely  many iterations of $\psi$ to a  transversal  $T$ in $\mathcal{S}_{\lambda}(Q_k)$, we will eventually  get a transversal $\Psi(T)$ in $\mathcal{S}_{\lambda}(P_k)$.  Thus, the map $\Psi$ is well-defined. 

\subsection{The correctness of the bijection}
The main objective of this subsection is to prove that the maps $\Phi$ and $\Psi$ are inverses of each other, thereby inducing  a one-to-one correspondence between $\mathcal{S}_{\lambda}(P_k)$ and $\mathcal{S}_{\lambda}(Q_k)$. 

\begin{lemma}\label{mainlem1}
	Let $k\geq 3$ and  $T\in \mathcal{S}_{\lambda}(P_k)$. For all $t\geq 1$, we have  $\psi(\phi^{t}(T))=\phi^{t-1}(T)$. 
\end{lemma}
\pf 
Suppose that at the $t$-th application of $\phi$ to $T$, the selected  submatrix   isomorphic to $Q_k$ is given by $G$ in   which  the squares containing $1's$ of   $G$  are labeled  by $b_1, b_2, \ldots, b_k$ from bottom to top,  and  the $1's$ of $G$ are labeled  by $a_1, a_2, \ldots, a_k$ from left to right. Clearly, $b_1$ is the square containing the lowest $1$ of $G$ and $b_k$ is square containing the topmost $1$ of $G$.    Recall that $\phi^{t}(T)$ is obtained from $\phi^{t-1}(T)$ by replacing $G$ with $\theta(G)$. Clearly, $\theta(G)$ is an occurrence of $P_k$.  Label the $1's$ of   $\theta(G)$ from left to right by $g_1, g_2, \ldots, g_k$.

 By the construction of $\theta(G)$, one can easily check that   $g_s$ is located at the same row as   $b_1$ for some $1\leq s\leq k-2$,   
$g_{k-1}$ is located at the same row as $b_k$ and the same column as $a_{k-1}$, and $g_k$ is located in the same column as $b_1$  as shown in  Figure \ref{fig:board1}.   Lemma \ref{Fact0} states that there are exactly $k-2$ $1's$ of $\phi^{t}(T)$ located in $F$.  It is apparent that such $k-2$ $1's$  coincide with    those $1's$  labeled by $g_1, g_2, \ldots, g_{s-1}, g_{s+1}, \ldots,  g_{k-1}$.  Again by  Lemma \ref{Fact0},   there are no $1's$ of  $\phi^{t}(T)$ located inside $E$.  
Assume that when we apply the map $\psi$ to $\phi^{t}(T)$, the selected submatrix is given by $H$.      Lemma~\ref{lemphi3} states that  $\phi^{t}(T)$ contains no $P_k$ that is higher than $\theta(G)$. Then    
by the rules specified in the selection    of $H$ in the definition of $\psi$, one can easily verify that $H=\theta(G)$.  

Recall that $\psi(\phi^{t}(T))$ is obtained from $\phi^{t}(T)$ by replacing $H=\theta(G)$ by $\theta'(\theta(G))$. 
In order to show that $\psi(\phi^{t}(T))=\phi^{t-1}(T)$, it remains to show that
$\theta'(\theta(G))=G$. By Lemma~\ref{lemphi4},  it follows that the  submatrix $G$ that we select in the application of $\phi$ to $\phi^{t-1}(T)$ is not isomorphic to $P_k$. This implies that $a_{k-1}$ is not positioned on the topmost row of $G$. Then  it follows from    Observation \ref{ob1} that   $\theta'(\theta(G))=G$ as desired, completing the proof. \qed

\begin{lemma}\label{mainlem2}
	Let $k\geq 3$ and  $T\in \mathcal{S}_{\lambda}(Q_k)$. Then  we have  $\phi(\psi^{t}(T))=\psi^{t-1}(T)$ for all $t\geq 1$. 
\end{lemma}
\pf  Assume that at the $t$-th  application of $\psi$ to $T$, the selected submatrix isomorphic to $P_k$ is given by $G$   whose squares containing  $1's$ are  labeled  by the same as those in the definition of $\psi$. To be more specific,  the leftmost  $k-1$ squares of $G$ containing $1's$ are labeled  by $b_1, b_2, \ldots, b_{k-1}$ from bottom to top and the rightmost square of $G$ containing a $1$ is labeled by $b_k$. 
Lemma \ref{lempsi4} tells us  that  $G$ is not isomorphic to $Q_k$. 
 For our convenience, we also label the $1's$ of $\theta'(G)$ by $g_1, g_2, \ldots, g_k$  from left to right.  Clearly, $\theta'(G)$ is a $Q_k$ and  $g_k$ is located at the square $x$ that is the lower right corner of $G$ as illustrated in Figure \ref{fig:board}.  
 Assume that at the application of $\phi$ to $\psi^t(T)$, the selected submatrix  isomorphic to $Q_k$ is given by $H$. Lemma \ref{lempsi3} ensures that $H$ is not lower than $\theta'(G)$.   
 
 Now we proceed to show that $H=\theta'(G)$.  Lemma \ref{Fact1} states that  there are exactly $k-1$ $1's$  of $\psi^{t}(T)$ that are positioned  in $F'$. Actually, such $1's$ coincide with those $1's$ labeled by  $g_1, g_2, \ldots, g_{k-1}$.  Again by Lemma \ref{Fact1}, 
   there are no $1's$ of    $\psi^{t}(T)$ that are located inside $E'$.  
  Recall that  $g_k$ is located at the square $x$.  Then by the rules specified in the selection of   $H$ in the defintion of the transformation $\phi$,   one can easily check that $H=\theta'(G)$ as desired.

Recall that $\phi(\psi^{t}(T))$ is obtained from $\psi^{t}(T)$ by replacing $H=\theta'(G)$ by $\theta (\theta'(G))$. 
In order to show that $\phi(\psi^{t}(T))=\psi^{t-1}(T)$, it remains to show that
$\theta(\theta'(G))=G$.      By  Observation \ref{ob2},  it remains to show that the rightmost $1$ of $G$ is not positioned at the lowest row of $G$. This  has been  justified by Lemma \ref{lempsi4}, completing the proof. \qed

\begin{theorem}
For $k\geq 3$, the maps $\Phi$ and $\Psi$ induce a bijection between $\mathcal{S}_{\lambda}(P_k)$ and $\mathcal{S}_{\lambda}(Q_k)$. 
\end{theorem}
\pf It remains to show that the map $\Phi$ and $\Psi$ are inverses of each other.   To this end, it suffices to show that  $\psi(\phi^t(T))=\phi^{t-1}(T)$ for any $t\geq 1$  and $T\in \mathcal{S}_{\lambda}(P_k)$ and  $\phi(\psi^t(T))=\psi^{t-1}(T)$ for any $t\geq 1$ and  $T\in \mathcal{S}_{\lambda}(Q_k)$. This has been  justified by Lemmas \ref{mainlem1} and  \ref{mainlem2}, completing the proof. \qed

\section*{Acknowledgments}
The authors thank the anonymous referees  for their insightful comments and suggestions.
 	The work  was supported by
	the National Natural
	Science Foundation of China (12471318).

\end{document}